\documentclass[12pt]{amsart}
\usepackage{epsfig,amsthm,amsmath,amsfonts,amssymb,color,nicefrac}
\usepackage[export]{adjustbox}
\newtheorem{theorem}{Theorem}

\numberwithin{equation}{section}

\theoremstyle{definition}

\theoremstyle{remark}

\newtheorem{assumption}[theorem]{Assumption}

\newcommand{\E}{\mathbb{E}}
\newcommand{\V}{\mathbb{V}}
\newcommand{\Cor}{\mathbb{C}\!\operatorname{or}}
\newcommand{\Cov}{\mathbb{C}\!\operatorname{ov}}
\renewcommand{\d}{\operatorname{d}\!}
\newcommand{\Rset}{\mathbb{R}}

\newcommand{\per}{{\rm per}}
\renewcommand{\Re}{\operatorname{Re}}

\newcommand{\red}[1]{\textcolor{black}{#1}}

\newcommand{\isdef}{\mathrel{\mathrel{\mathop:}=}}

\begin{document}

\title[Rapid computation of far-field statistics]
{Rapid computation of far-field statistics for random obstacle scattering}
\author{Helmut Harbrecht}
\address{Helmut Harbrecht, 
Departement Mathematik und Informatik, 
Universit\"at Basel, 
Spiegelgasse 1, 4051 Basel, Schweiz.}
\email{helmut.harbrecht@unibas.ch}
\author{Nikola Ili\'c}
\address{Nikola Ili\'c,
Departement Mathematik und Informatik, 
Universit\"at Basel, 
Spiegelgasse 1, 4051 Basel, Schweiz.}
\email{n.ilic@unibas.ch}
\author{Michael D.~Multerer}
\address{
Michael D.~Multerer,
Institute of Computational Science,
USI Lugano,
Via Giuseppe Buffi 13, 6900 Lugano, Schweiz}
\email{michael.multerer@usi.ch}
\date{}
\begin{abstract}
In this article, we 
consider the numerical approximation of far-field statistics for
acoustic scattering problems in the case of random obstacles.  
In particular, we consider the computation 
of the expected far-field pattern and the expected scattered 
wave away from the scatterer as well as the computation of the
corresponding variances.
To that end, we introduce an artificial interface, which almost
surely contains all realizations of the random scatterer.
At this interface, we directly approximate the second order 
statistics, i.e., the expectation and the variance, of the Cauchy data
by means of boundary integral equations. From these quantities, we
are able to rapidly evaluate statistics of the scattered wave 
everywhere in the exterior domain, including the 
expectation and the variance of the far-field. By employing a low-rank
approximation of the Cauchy data's two-point correlation function,
we drastically reduce the cost of 
the computation of the scattered wave's variance. 
Numerical results are provided in order to demonstrate 
the feasibility of the proposed approach.
\end{abstract}
\maketitle 

\section{Introduction}\label{sec:intro}
The propagation of an acoustic wave in a homogeneous, 
isotropic, and inviscid fluid is approximately described by a 
velocity potential $U({\bf x},t)$ satisfying the wave equation
\[
U_{tt} =  c^2\Delta U.
\]
Here, $c$ denotes the speed of sound, ${\bf v} = \nabla U$ is the
velocity field, and $p = -U_t$ is the pressure, see \cite[Chapter 3]{CK2} for
instance. If $U$ is time harmonic, that is 
\[
  U({\bf x},t) = \Re\big(u({\bf x})e^{-i\omega t}\big), \quad\omega > 0, 
\]  
in complex notation, then the complex-valued space-dependent 
function $u$ satisfies the Helmholtz equation
\[
  \Delta u + \kappa^2 u = 0\quad\text{in $\mathbb{R}^d\setminus\overline{D}$},
\]
where $D\subset\Rset^d$, \(d=2,3\), corresponds to an obstacle and 
$\kappa = \omega/c$ is the wavenumber. We assume that 
$D$ is a bounded and simply connected domain, having a smooth boundary 
$\Gamma = \partial D$. For \emph{sound-soft} obstacles the pressure $p$ 
vanishes on $\Gamma$, which leads to the Dirichlet boundary condition
\[
  u=0\quad\text{on $\Gamma$}.
\]
We shall consider the situation that the total wave 
\[u = u_{\operatorname{inc}}+u_{\operatorname{s}}\] is comprised of a
known incident plane wave $u_{\operatorname{inc}}({\bf x}) = e^{i\kappa \langle{\bf d},{\bf x}\rangle}$ 
with direction ${\bf d}\in\mathbb{R}^d$, where $\|{\bf d}\|_2=1$, and the scattered wave $u_{\operatorname{s}}$. 
Then, if we impose the \emph{Sommerfeld radiation condition}
\begin{equation}	\label{eq:pde3}
  \sqrt{r}\bigg(\frac{\partial u_{\operatorname{s}}}{\partial r}-i\kappa u_{\operatorname{s}}\bigg)
  	\to 0\ \text{as}\ r \isdef \|{\bf x}\|_2\to\infty,
\end{equation}
for the scattered wave, we obtain a unique solution to the acoustic scattering problem
\begin{equation}\label{eq:pde1}
 \begin{aligned}
  \Delta u + \kappa u= 0\quad&\text{in}\ \mathbb{R}^d\setminus\overline{D},\\
  	u= 0\quad&\text{on}\ \Gamma,\\
	  \sqrt{r}\bigg(\frac{\partial u_{\operatorname{s}}}{\partial r}-i\kappa u_{\operatorname{s}}\bigg)
  	\to 0\quad&\text{as}\ r = \|{\bf x}\|_2\to\infty,
  \end{aligned}
\end{equation}
see \cite[Chapter 3]{CK2}.
In particular, the Sommerfeld radiation condition implies the asymptotic behavior
\begin{equation}	\label{eq:pde4}
  u_{\operatorname{s}}({\bf x}) = \frac{e^{i\kappa r}}{r} 
  	\left\{u_\infty\bigg(\frac{\bf x}{r}\bigg)
	+ \mathcal{O}\bigg(\frac{1}{r}\bigg)\right\},
		\quad r\to\infty.
\end{equation}
Herein, the function 
\[
u_\infty\colon\mathbb{S}^1\isdef\{\hat{\bf x}\in\mathbb{R}^d:\|\hat{\bf x}\|_2=1\}
\to\mathbb{C}
\] 
is called the \emph{far-field pattern}, which is always analytic in accordance with \cite[Chapter 6]{CK2}.

In this article, we consider the situation that the scatterer $D$ is 
randomly shaped, i.e., $D=D({\bf y})$ for a random parameter \({\bf y}\in\square\isdef[-1,1]^\mathbb{N}\). 
Hence, the scattered  wave itself becomes a random field $u_{\operatorname{s}}({\bf y})$. We will model
a class of random domains and compute the associated expected 
scattered wave $\E[u_{\operatorname{s}}]$ and also the expected far-field $\E[u_\infty]$.
\red{Instead of employing the domain mapping method, which maps
the deformed scatterer onto a fixed reference domain, 
as in e.g.\ \cite{CNT16,HPS16,H3S,JSZ17,xiu}, or a fictitious domain approach as in
\cite{CK07}, we will compute 
all samples for the deformed scatterer by means of the boundary element
method.} 
This approach is much cheaper since we do not require a very fine triangulation 
for \(D\) in order to ensure that the domain deformation 
field is properly resolved. Consequently, 
we are also able to deal with large variations without the need of
a very fine discretization.

Furthermore, we derive a means to compute the scattered wave's second 
order statistics in a deterministic fashion from its Cauchy data's second 
order statistics on an artificial, deterministic interface $\Sigma$, which 
almost surely contains the domain $D({\bf y})$.
By the application of a low-rank approximation for the correlation 
function, we are able to considerably decrease the cost for the 
computation of the expected scattered 
field and its variance. The advantages of the proposed approach are 
thus as follows: 
\begin{itemize}
\item[(i)] The use of boundary integral equations facilitates a 
straightforward treatment of the unbounded exterior domain. Especially, it 
avoids expensive mesh generation procedures in case of strongly 
varying scatterers. 
\item[(ii)] Since the artificial interface is bounded and 
has one dimension less than to the exterior domain, the impact 
of the high dimensionality of the random scattering problem is
drastically reduced.
\end{itemize}

We like to emphasize that the present approach is also suitable 
to treat sound-hard scatterers, where the Dirichlet boundary condition in \eqref{eq:pde1}
becomes a Neumann condition. In addition, scatterers with a different diffractive 
index $\kappa$ can be considered. The latter leads to a transmission condition
at the scatterer's surface instead of a boundary condition. The presented ideas remain valid
in this situation except for modifying the boundary integral equations accordingly. \red{
Moreover, although we focus on \(d=2\) in the numerical examples, 
all concepts can be transferred to \(d=3\) in a straightforward manner. However,
technicalities will increase.}

The rest of the article is organized as follows. In Section \ref{sec:bie},
we introduce the formulation of the scattering problem under consideration 
in the case of a deterministic scatterer by means of boundary integral equations.
In particular, we provide a representation of the total wave and the 
far-field pattern. Then, in Section \ref{sec:random},
we consider a representation of random scatterers in terms of random vector
fields. Moreover, we provide an explicit description for the two dimensional
situation, which is used later on in the numerical examples.
Section \ref{sec:randScat} deals with the random scattering problem. 
Here, we derive expressions for 
the scattered wave's expectation and variance, including the far-field pattern.
Section~\ref{sec:results} is dedicated to numerical results which quantify and
qualify our approach. The boundary integral equations are discretized by 
the Nystr\"om method which converges exponentially in case of analytic
boundaries. Especially, we discuss the efficient computation of the 
scattered wave's variance by using a low-rank approximation. Finally, 
in Section~\ref{sec:conclusion}, we state concluding remarks.

\section{Boundary integral representation of the scattering problem}\label{sec:bie}
\subsection{Computing the scattered wave}
We shall recall the solution of the boundary value problem 
\eqref{eq:pde1} by means of boundary integral
equations. For the sake of simplicity in representation, 
we assume here that the domain $D$ is fixed with a smooth
boundary $\Gamma=\partial D$. 

We introduce the acoustic single layer operator
\[
  \mathcal{V}\colon H^{-\nicefrac{1}{2}}(\Gamma)\to H^{\nicefrac{1}{2}}(\Gamma),\quad
  \mathcal{V}\rho\isdef\int_\Gamma 
  	\Phi(\cdot,{\bf z})\rho({\bf z})\d\sigma_{\bf z}
\]
and the acoustic double layer operator
\[
  \mathcal{K}\colon L^2(\Gamma)\to L^2(\Gamma),\quad
  \mathcal{K}\rho\isdef\int_\Gamma 
  	\frac{\partial\Phi(\cdot,{\bf z})}{\partial{\bf n}_{\bf z}}\rho({\bf z})\d\sigma_{\bf z}.
\]
Herein, $\Phi(\cdot,\cdot)$ denotes the fundamental 
solution of the Helmholtz equation. It is given by
\[
  \Phi({\bf x}, {\bf x}') = \begin{cases}
  \displaystyle{\frac{i}{4}H_0^{(1)}(\kappa\|{\bf x}-{\bf x}'\|_2)}, & d = 2,\\[1em]
  \displaystyle{\frac{e^{i\kappa\|{\bf x}-{\bf x}'\|_2}}{4\pi\|{\bf x}-{\bf x}'\|_2}}, & d=3.
  \end{cases}
\]
where $H_0^{(1)}$ denotes the zeroth order Hankel function 
of the first kind.

Then, if the incident wave is given by 
\[u_{\operatorname{inc}}({\bf x}) = e^{i\kappa\langle{\bf d},
{\bf x}\rangle}\] for some direction \({\bf d}\in\mathbb{R}^d\), 
the Neumann data of the total wave $u=u_{\operatorname{inc}}+u_{\operatorname{s}}$ at the 
boundary $\Gamma$ can be determined by the boundary integral equation
\begin{equation}\label{eq:BIE1}
  \left(\frac{1}{2} + \mathcal{K^\star} - i\eta\mathcal{V}\right)
  	\frac{\partial u}{\partial {\bf n}} = \frac{\partial u_{\operatorname{inc}}}{\partial {\bf n}} - i\eta u_{\operatorname{inc}}
		\quad\text{on $\Gamma$},
\end{equation}
where \(\mathcal{K}^\star\) denotes the adjoint double-layer operator, \({\bf n}\) the
outward pointing normal vector and \(\eta\in\mathbb{R}\setminus\{0\}\)
is chosen such that \(\eta\operatorname{Re}(\kappa)>0\),
see \cite{BW,CK2}. 

From the Cauchy data of $u$ at $\Gamma$, we 
can determine the scattered wave $u_{\operatorname{s}}$ in any point in the
exterior of the scatterer by applying the potential evaluation
\begin{equation}\label{eq:solution1}
  u_{\operatorname{s}}({\bf x}) = \int_\Gamma\Phi({\bf x}, {\bf z})
	\frac{\partial u}{\partial{\bf n}}({\bf z})\d\sigma_{\bf z},
	\quad {\bf x}\in\mathbb{R}^d\setminus\overline{D}.
\end{equation}
By letting $\|{\bf x}\|_2$ tend to infinity in \eqref{eq:solution1}, 
we derive a closed expression for the far-field of the total 
wave $u$. Namely, the far-field at a point $\widehat{\bf x}\in
\mathbb{S}^1$ is given in accordance with
\begin{equation}\label{eq:far-field1}
  u_\infty(\hat{\bf x}) = \int_\Gamma\Phi_\infty(\hat{\bf x},{\bf z})
	\frac{\partial u}{\partial{\bf n}}({\bf z})\d\sigma_{\bf z}.
\end{equation}
Herein, the far-field kernel $\Phi_\infty(\cdot,\cdot)$ is given according to
\[
  \Phi_\infty(\hat{\bf x},{\bf z}) = 
  \begin{cases} \displaystyle{\frac{e^{i\pi /4}}{\sqrt{8\kappa\pi}}
  	e^{-i\kappa\langle\hat{\bf x},{\bf z}\rangle}}, & d=2,\\[1em]
	\displaystyle{\frac{1}{4\pi}e^{-i\kappa\langle\hat{\bf x},{\bf z}\rangle}},& d=3.
	\end{cases}
\]

\subsection{Alternative representation of the scattered wave}
We shall introduce the sphere
\[
\Sigma \isdef\{{\bf x}\in\mathbb{R}^d:\|{\bf x}\|_2=R\}
\]
of radius $R>0$, being sufficiently large to guarantee that $\Sigma$ encloses 
the domain $D$. By differentiating \eqref{eq:solution1}, it is seen that 
the gradient of the scattered wave can simply be computed by
\[
  \nabla u_{\operatorname{s}}({\bf x}) = \int_\Gamma \nabla_{\bf x} \Phi({\bf x}, {\bf z})
	\frac{\partial u}{\partial{\bf n}}({\bf z})\d\sigma_{\bf z},
	\quad {\bf x}\in\Sigma.
\]
Thus, we can compute the Cauchy data of the scattered 
wave at the artificial interface $\Sigma$. Especially, it holds
\[
  \frac{\partial u_{\operatorname{s}}}{\partial {\bf n}}({\bf x}) = \int_\Gamma \frac{\partial\Phi({\bf x}, {\bf z})}
  	{\partial {\bf n}_{\bf x}}\frac{\partial u}{\partial{\bf n}}({\bf z})\d\sigma_{\bf z},
	\quad {\bf x}\in\Sigma,
\]
where ${\bf n}_{\bf x} = \nicefrac{\bf x}{\|{\bf x}\|_2}$ is the 
outward pointing normal at ${\bf x}\in\Sigma$. 

For any ${\bf x}\in\mathbb{R}^d$ with $\|{\bf x}\|_2>R$, 
we can now either use the representation formula \eqref{eq:solution1} 
or the representation formula
\begin{equation}\label{eq:solution2}
  u_{\operatorname{s}}({\bf x}) = \int_\Sigma \bigg\{\Phi({\bf x}, {\bf z})
	\frac{\partial u_{\operatorname{s}}}{\partial{\bf n}}({\bf z})+
	\frac{\partial\Phi({\bf x}, {\bf z})}{\partial {\bf n}_{\bf z}}
		u_{\operatorname{s}}({\bf z})\bigg\}\d\sigma_{\bf z}
\end{equation}
to evaluate the scattered wave $u_{\operatorname{s}}$ at any point ${\bf x}
\in\mathbb{R}^d$ with $\|{\bf x}\|_2> R$. In particular, letting \(R\to\infty\), we obtain for the far-field the formula
\begin{equation}\label{eq:far-field2}
  u_\infty(\hat{\bf x}) = \int_\Sigma\bigg\{\Phi_\infty(\hat{\bf x},{\bf z})
	\frac{\partial u}{\partial{\bf n}}({\bf z}) 
	+ \frac{\partial\Phi_\infty({\bf x}, {\bf z})}{\partial {\bf n}_{\bf z}}
		u({\bf z})\bigg\}\d\sigma_{\bf z},\quad\hat{\bf x}\in\mathbb{S}^1.
\end{equation}
As we will see, the major advantage of \eqref{eq:solution2} and 
\eqref{eq:far-field2} over \eqref{eq:solution1} and \eqref{eq:far-field1} in case
of a random scatterer
is that the sphere $\Sigma$ has a fixed shape. 

We remark that an artificial interface being different from a circle 
can of course be chosen as well. 
\section{Random obstacles}\label{sec:random}
\subsection{Representation of random domains}
\red{In this section, we introduce a description of random obstacles by
means of random vector fields, as they have originally been considered
in \cite{HPS16} in the context of the domain mapping method.
To that end, let \((\Omega,\mathcal{A},\mathbb{P})\) denote
a complete and separable probability space with \(\sigma\)-algebra
\(\mathcal{A}\) and probability measure \(\mathbb{P}\). Here, complete
means that \(\mathcal{A}\) contains all \(\mathbb{P}\)-null sets.
For a given real or complex Banach space $\mathcal{X}$, we introduce the 
\emph{Lebesgue-Bochner space}
$L_{\mathbb{P}}^p(\Omega;\mathcal{X})$, 
$1\le p\le\infty$, which consists of
all equivalence classes of strongly measurable functions 
$v\colon\Omega\to\mathcal{X}$ with bounded norm
\[
  \|v\|_{L_{\mathbb{P}}^p(\Omega;\mathcal{X})}\isdef
  \begin{cases}
    \displaystyle{
      \left(\int_\Omega \|v(\cdot,\omega)\|_{\mathcal{X}}^p \,\d\mathbb{P}(\omega)\right)^{1/p}
    },
    & p<\infty \\[2ex]
    \displaystyle{
      \operatorname*{ess\,sup}_{\omega\in\Omega} \|v(\cdot,\omega)\|_{\mathcal{X}}
    },
    & p=\infty,
  \end{cases}
\]
If $p=2$ and $\mathcal{X}$ is a separable Hilbert space, then the
Lebesgue-Bochner space $L_{\mathbb{P}}^p(\Omega;\mathcal{X})$ is isomorphic to the tensor product space
$L_{\mathbb{P}}^2(\Omega)\otimes \mathcal{X}$.
For more details on Lebesgue-Bochner spaces, we refer the reader to \cite{HP57}.}

\red{
For \(p\geq 2\) and a given \emph{random field} \(v\in L^p(\Omega;\mathcal{X})\), 
we can introduce the \emph{expectation}
\[
 \E[v]({\bf x})\isdef\int_\Omega v({\bf x},\omega)\d\mathbb{P}(\omega)
\]
and the \emph{variance}
\[
\V[v]({\bf x})\isdef\int_\Omega v({\bf x},\omega)\overline{v({\bf x},\omega)}\d\mathbb{P}(\omega)-\E[v]({\bf x})\overline{\E[v]({\bf x})}.
\]
With straighforward modifications, these definitions remain valid for real valued random fields.}

\red{
Now, to define the random obstacle \(D(\omega)\subset\mathbb{R}^d\), 
we assume the existence of a nominal obstacle \(D_0\subset\mathbb{R}^d\), 
with boundary \(\Gamma_0\isdef\partial D_0\), and of a 
\emph{uniform \(C^1\)-diffeomorphism}
\[{\bf V}\colon\overline{D_0}\times\Omega\to\mathbb{R}^d,
\]
i.e.\ there holds
\begin{equation}\label{eq:unifVfield}
\|{\bf V}(\omega)\|_{C^1(\overline{D_0};\mathbb{R}^d)},\|{\bf V}^{-1}(\omega)\|_{C^1(\overline{D_0};\mathbb{R}^d)}\leq C_{\operatorname{uni}}\quad\text{for \(\mathbb{P}\)-a.e.\ }\omega\in\Omega,
\end{equation}
such that \(D(\omega)\) is implicitly given by the relation
\[
D(\omega)={\bf V}(D_0,\omega).
\]
Consequently, we obtain
\[
\Gamma(\omega)\isdef\partial D(\omega)={\bf V}(\Gamma_0,\omega).
\]
Moreover, we define the \emph{hold-all domain} \(\mathcal{D}\) 
according to
\begin{equation}\label{eq:holdall}
\mathcal{D}\isdef\bigcup_{\omega\in\Omega}D(\omega).
\end{equation}
Due to \eqref{eq:unifVfield}, it holds \({\bf V}\in 
L^\infty\big(\Omega;C^1(\overline{D_0})\big)\subset  
L^2\big(\Omega;C^1(\overline{D_0})\big)\). Hence, 
the vector field \({\bf V}\) can be represented by 
a \emph{Karhunen-Lo\`eve expansion}, 
cf.\ \cite{Loeve}, of the form
 \begin{equation}\label{eq:randomVectorField}
{\bf V}({\bf x},\omega)=\E[{\bf V}]({\bf x})+
\sum_{k=1}^\infty{\bf V}_k({\bf x}){Y}_k(\omega).
\end{equation}
An efficient way to compute the Karhunen-Lo\`eve 
expansion if the mean and the covariance function of 
the random (vector) field under consideration are known 
is given by the \emph{pivoted Cholesky decomposition}. 
This accounts particularly for random vector fields, 
see \cite{HPS,HPS15,Mul18}.}

\red{
The anisotropy which is induced by the spatial parts \(\{{\bf V}_k\}_k\), 
describing the fluctuations around the nominal value $\mathbb{E}[{\bf V}]({\bf x})$ 
in \eqref{eq:randomVectorField}, is encoded by the quantities
\begin{equation}\label{eq:GammaK}
\gamma_k\isdef\|{\bf V}_k\|_{W^{1,\infty}(D_0;\mathbb{R}^d)}.
\end{equation}
For our modeling, we shall also make the following common assumptions.
\begin{assumption}\label{ass:KL}\
    \begin{itemize}
    \item[(i)] The random variables \(\{Y_k\}_k\) take values in \([-1,1]\).
    \item[(ii)] The random variables \(\{Y_k\}_k\) are independent and uniformly distributed,
    i.e.\ \(Y_k\sim\mathcal{U}(-1,1)\).
    \item[(iii)] The sequence \(\{{\gamma}_k\}_k\) is at least in \(\ell^1(\mathbb{N})\).
\end{itemize}
\end{assumption}
We remark that it holds without loss of generality \(\E[{\bf V}]({\bf x})={\bf x}\), otherwise
we have to choose an appropriate reparametrization.
Moreover, identifying the random variables by their image \({\bf y}\in\square\isdef[-1,1]^{\mathbb{N}}\),
we end up with the representation
\begin{equation}\label{eq:ParVfield}
    {\bf V}({\bf x},{\bf y})={\bf x}+
    \sum_{k=1}^\infty{\bf V}_k({\bf x}){y}_k.
\end{equation}
The corresponding image measure \(\mu\) is given by the product of the push forward measure 
\(\nu=\nicefrac{\d y}{2}\) according to \(\mu\isdef\otimes_{k=1}^\infty\nu\).
The Jacobian of \({\bf V}\) with respect to the spatial variable ${\bf x}$ is thus given by
\[
{\bf J}({\bf x},{\bf y})={\bf I}+
\sum_{k=1}^\infty{\bf V}_k'({\bf x}){y}_k.
\]
}

\red{
The uniformity condition implies that the functional determinant \(\det {\bf J}({\bf x},{\bf y})\)
is either uniformly positive or negative, see \cite{HPS16} for more details. Hence, we may assume
without loss of generality
\[
0<c\leq\det {\bf J}({\bf x},{\bf y})\leq C<\infty
\]
for every \({\bf x}\in D_0\) and almost every \({\bf y}\in\square\), 
where \(c,C>0\) are some constants. For a sufficiently fine discretization of the domain \(D_0\),
this property carries over to the finite element approximation of \({\bf V}\) and
\({\bf J}\), respectively. Consequently, a
quasi-uniform mesh will always be mapped to a quasi-uniform mesh by \({\bf V}\). Again, we refer to \cite{HPS16} for the details.
Hence, under the uniformity condition \eqref{eq:unifVfield} no remeshing of \(D({\bf y})\)
or \(\Gamma({\bf y})\) is necessary for different realizations of \({\bf y}\in\square\). However, we emphasize that the discretization of \({\bf V}\)
and hence the mesh for \(D_0\) needs to be sufficiently fine in order to guarantee 
\eqref{eq:unifVfield} also for the discretized  random vector field.}

\subsection{Star-like obstacles in two spatial dimensions}
\label{subseq:starlike}
In the numerical examples, we will restrict ourselves to star-like
scatterers in two spatial dimensions.
Especially, since we consider a boundary integral 
approach to solve the Helmholtz equation, we shall only define here
the random vector field only for the boundary \(\Gamma_0\).

Without loss of generality, we assume that the scatterer \(D(\omega)\)
is  star-like with respect to the origin \({\bf 0}\in\mathbb{R}^2\). 
Then, we  can represent its boundary $\Gamma({\omega})$ 
by a parametrization of the form
\begin{equation}\label{eq:boundrep}
  \boldsymbol\gamma\colon[0,2\pi]\times\Omega\to\mathbb{R}^2,\quad
 \boldsymbol\gamma(\phi)=r(\phi,\omega){\bf e}_r(\phi),
\end{equation}
where \[{\bf e}_r(\phi) = \begin{bmatrix}\cos(\phi)\\ 
\sin(\phi)\end{bmatrix}\]
denotes the radial 
direction and 
radius function \(r(\phi,\omega)\) is a real valued random 
field \(r\in L^2\big(\Omega;C^n_{\text{per}}([0,2\pi])\big)\).
As before, the Karhunen-Lo\`eve expansion of \(r\) can be computed
if \(\E[r](\phi)\) and 
\[
\Cov[r]({\phi},{\phi}')\isdef
\int_\Omega\big(r(\phi,\omega)-\E[r](\phi)\big)
{\big(r(\phi',\omega)-\E[r](\phi')\big)}\d\mathbb{P}(\omega).
\]
are known. However, we will assume here that \(r\) is of the particular form
\begin{equation}\label{eq:RandomRadius}
  r(\phi,{\bf y}) = r_0(\phi) + \sum_{k=1}^\infty \big\{a_{2k-1} y_{2k-1} \sin(k\phi)	
  	+ a_{2k} y_{2k} \cos(k\phi)\big\},\quad{\bf y}\in\square.
\end{equation}
In this case, the spatial regularity is entirely encoded in the coefficients
\(a_k\), \(k\geq 0\), see \cite{H3S}.
By construction, the random fluctuations of the radius 
\eqref{eq:RandomRadius} are centered, i.e., their mean 
vanishes, and we conclude
\[
  \mathbb{E}[r](\phi) = \int_{\square} r(\phi,{\bf y}) \d\mu = r_0(\phi).
\]

In order to guarantee that each realization ${\bf y}\in [-1,1]^\mathbb{N}$ results in a boundary of a valid domain \(D({\bf y})\), we shall further assume that
\begin{equation*}
  0<\underline{r}\le r(\phi,{\bf y})\le\overline{r}<\infty
  	\quad\text{for all $\phi\in[0,2\pi]$ and ${\bf y}\in\square$}.
\end{equation*}
Moreover, it is assumed that $r_0\in C_{\per}^2([0,2\pi])$ as well as 
that the sequence $(a_k)_k$ decays sufficiently fast to ensure $r(\cdot,{\bf y})
\in C_{\per}^2([0,2\pi])$ for all ${\bf y}\in\square$. \red{This for example 
guaranteed if \(|a_k|\leq ck^{-(3+\varepsilon)}\) for any \(\varepsilon>0\) 
and some constant \(c>0\).}

The random boundary $\Gamma({\bf y})$ is hence given by
\[
  \Gamma({\bf y}) = \big\{{\boldsymbol\gamma}(\phi,{\bf y})=r(\phi,{\bf y}){\bf e}_r(\phi)\in\mathbb{R}^2:
  	\phi\in[0,2\pi]\big\}.
\]
Consequently, there holds \({\bf V}|_{\Gamma_0}={\boldsymbol\gamma}\), where
\(\Gamma_0=r_0{\bf e}_r\). In order to provide a description of the 
random vector field for the interior of the scatterer, 
which is for example needed in a domain mapping approach, 
a suitable extension of \({\bf V}\) to \(D_0\) has to be defined as in 
e.g.\ \cite{HPS16,H3S}. 
However, we remark that in the approach presented here, no knowledge of 
\({\bf V}\) inside the scatterer is required.

\section{The random scattering problem}\label{sec:randScat}
\subsection{Problem formulation}
Now, having a suitable description of the random scatterer \(D({\bf y})\)
at our disposal, we can define the random scattering problem under consideration
Let \(u_{\operatorname{inc}}\) denote the incident wave. Then, the boundary 
value problem for the total field $u({\bf y}) = u_{\operatorname{s}}({\bf y})+u_{\operatorname{inc}}$
for a given ${\bf y}\in\square$
reads 
\begin{equation}\label{eq:spde}
 \begin{aligned}
  \Delta u({\bf y}) + \kappa u({\bf y}) = 0\quad&\text{in}\ \mathbb{R}^d\setminus\overline{D({\bf y})},\\
  	u({\bf y}) = 0\quad&\text{on}\ \Gamma({\bf y}),\\
	  \sqrt{r}\bigg(\frac{\partial u_{\operatorname{s}}}{\partial r}-i\kappa u_{\operatorname{s}}\bigg)
  	\to 0\quad&\text{as}\ r = \|{\bf x}\|_2\to\infty.
  \end{aligned}
\end{equation}
By the construction of \(\Gamma({\bf y})\), the random scattering problem \eqref{eq:spde} exhibits
a unique solution for each realization \({\bf y}\in\square\) of the random parameter.
\red{Moreover, it has been shown in \cite{H3S} for the case of the Helmholtz transmission problem in two spatial dimensions that the
total wave \(u({\bf y})\) in the large wavelength regime exhibits 
an analytic extension with respect to the parameter \({\bf y}\in\square\) 
into a certain region of the complex plane. Hence, given that \(\kappa\) 
is sufficiently small, we may employ higher order quadrature methods, 
like the quasi-Monte Carlo methods, see e.g.\ \cite{Caf98}, or sparse 
quadrature methods, see e.g. \cite{HHPS18,H3S}, for the computation 
of quantities of interest.}


\subsection{Expected scattered wave}
We can compute the scattered wave's expectation for 
a given point ${\bf x}\in\mathbb{R}^d$ via the potential 
evaluation \eqref{eq:solution1}, which leads to
\begin{equation}\label{eq:expectation1}
  \E[u_{\operatorname{s}}]({\bf x}) = \E\bigg[\int_{\Gamma({\bf y})}\Phi({\bf x}, {\bf z})
	\frac{\partial u_{\operatorname{s}}}{\partial{\bf n}}({\bf z},\cdot)\d\sigma_{\bf z}\bigg].
\end{equation}
Of course, \eqref{eq:expectation1} makes only sense
if ${\bf x}\not\in\mathcal{D}$ since otherwise there might be
instances ${\bf y}\in\square$ 
such that ${\bf x}\in D({\bf y})$, i.e., the point ${\bf x}$ does 
not lie outside the scatterer almost 
surely, compare \eqref{eq:holdall}.
In what follows, we assume that \(R>0\) is chosen such that
\[
\mathcal{D}\subset B_R({\bf 0})\isdef\{{\bf x}\in\mathbb{R}^d: \|x\|_2<R\}\quad\text{and set }
\Sigma=\partial B_R({\bf 0}).
\] 

Then, if we want to compute the expectation
in many points, it is much more efficient to exploit the
artificial but fixed boundary $\Sigma$ and to consider 
an expression similar to \eqref{eq:solution2}. 
\red{However, we remark that an application 
of the mapping approach as in \cite{H3S} is 
better suited if statistics of the scattered wave in 
the vicinity of the scatterer are of interest.}

For any ${\bf x}
\in\mathbb{R}^d$ with $\|{\bf x}\|_2>R$, it holds
\begin{equation}\label{eq:expectation2}
  \E[u_{\operatorname{s}}]({\bf x}) = \int_\Sigma \bigg\{\Phi({\bf x}, {\bf z})
	\E\bigg[\frac{\partial u_{\operatorname{s}}}{\partial{\bf n}}\bigg]({\bf z})+
	\frac{\partial\Phi({\bf x}, {\bf z})}{\partial {\bf n}_{\bf z}}
		\E[u_{\operatorname{s}}]({\bf z})\bigg\}\d\sigma_{\bf z}.
\end{equation}
Therefore, the scattered wave's expectation is completely
encoded in the Cauchy data at the artificial boundary $\Sigma$.
This means that we only need to compute the expected Cauchy 
data 
\begin{equation}\label{eq:expectation3}
  \mathbb{E}[u_{\operatorname{s}}] = \int_{\square} 
  	\bigg\{\int_{\Gamma({\bf y})}\Phi({\bf x}, {\bf z})
	\frac{\partial u}{\partial{\bf n}}({\bf z},{\bf y})\d\sigma_{\bf z}\bigg\}\d\mu
\end{equation}
and
\begin{equation}\label{eq:expectation4}
  \mathbb{E}\bigg[\frac{\partial u_{\operatorname{s}}}{\partial {\bf n}}\bigg]
  	= \int_{\square} 
  		\bigg\{\int_{\Gamma({\bf y})}\frac{\partial\Phi({\bf x}, {\bf z})}{\partial {\bf n}_{\bf z}}
			u({\bf z},{\bf y})\d\sigma_{\bf z}\bigg\}\d\mu
\end{equation}
of the scattered wave at the artificial boundary $\Sigma$. 

In complete analogy to \eqref{eq:expectation2}, the expected 
far-field pattern is likewise computed by using \eqref{eq:far-field2}:
\[
  \E[u_\infty](\hat{\bf x}) = \int_\Sigma\bigg\{\Phi_\infty(\hat{\bf x},{\bf z})
	\E\bigg[\frac{\partial u_{\operatorname{s}}}{\partial{\bf n}}({\bf z})\bigg]
	+ \frac{\partial\Phi_\infty({\bf x}, {\bf z})}{\partial {\bf n}_{\bf z}}
		\E[u_{\operatorname{s}}]({\bf z})\bigg\}\d\sigma_{\bf z}.
\]

\subsection{Computing the solution's variance}
The variance $\V[u_{\operatorname{s}}]$ of the scattered wave $u_{\operatorname{s}}$ at 
a point ${\bf x}\not\in B_R({\bf 0})$ depends nonlinearly on 
the Cauchy data of $u_{\operatorname{s}}$ at the artificial interface $\Sigma$. 
Nonetheless, we can employ the fact that the variance is the trace 
of the covariance function:
\begin{equation}\label{eq:var}
  \V[u_{\operatorname{s}}]({\bf x}) = \Cov[u_{\operatorname{s}}]({\bf x},{\bf x}')\big|_{{\bf x}={\bf x}'}\\
  	= \Cor[u_{\operatorname{s}}]({\bf x},{\bf x}')\big|_{{\bf x}={\bf x}'} - |\E[u_{\operatorname{s}}]({\bf x})|^2.
\end{equation}
The covariance function is given by
\begin{align*}
  \Cov[u_{\operatorname{s}}]({\bf x},{\bf x}') &= \E\Big[\big(u_{\operatorname{s}}({\bf x},\cdot)-\E[u_{\operatorname{s}}]({\bf x})\big)
  	\overline{\big(u_{\operatorname{s}}({\bf x}',\cdot)-\E[u_{\operatorname{s}}]({\bf x}')\big)}\Big]\\
	&=\E\big[u_{\operatorname{s}}({\bf x},\cdot)\overline{u_{\operatorname{s}}({\bf x'},\cdot)}\big]-\E[u_{\operatorname{s}}]({\bf x})\overline{\E[u_{\operatorname{s}}]({\bf x}')},
\end{align*}
and, hence, 
\[\Cor[u_{\operatorname{s}}]({\bf x},{\bf x}') = \E\big[u_{\operatorname{s}}({\bf x},\cdot)\overline{u_{\operatorname{s}}({\bf x'},\cdot)}\big].
\]
Hence, the two-point correlation function is 
a higher-dimensional object. Fortunately, it depends only linearly on 
the second moments of the Cauchy data of the scattered wave 
on the artificial interface $\Sigma$, which greatly simplifies its
computation. Namely, defining for ${\bf x},{\bf x}'\in\Sigma$ the
quantities
\begin{align*}
  &\Cor[u_{\operatorname{s}}]({\bf x},{\bf x}')\\
  &\qquad= \E\bigg[\bigg(\int_{\Gamma({\bf y})}\!\!\!\!\!\Phi({\bf x}, {\bf z})
	\frac{\partial u_{\operatorname{s}}}{\partial{\bf n}}({\bf z},{\bf y})\d\sigma_{\bf z}\bigg)
	\overline{\bigg(\int_{\Gamma({\bf y})}\!\!\!\!\!\Phi({\bf x}', {\bf z})
	\frac{\partial u_{\operatorname{s}}}{\partial{\bf n}}({\bf z},{\bf y})\d\sigma_{{\bf z}}\bigg)}\bigg],\\
 & \Cor\bigg[\frac{\partial u_{\operatorname{s}}}{\partial{\bf n}}\bigg]({\bf x},{\bf x}')\\
  &\qquad
  = \E\bigg[\bigg(\int_{\Gamma({\bf y})}\!\!\!\!\!\frac{\partial\Phi({\bf x}, {\bf z})}{\partial {\bf n}_{\bf z}}
		u_{\operatorname{s}}({\bf z},{\bf y})\d\sigma_{\bf z}\bigg)
	\overline{\bigg(\int_{\Gamma({\bf y})}\!\!\!\!\!\frac{\partial\Phi({\bf x}', {\bf z})}{\partial {\bf n}_{\bf z}}
		u_{\operatorname{s}}({\bf z},{\bf y})\d\sigma_{{\bf z}}\bigg)}\bigg],
\end{align*}
and
\begin{align*}
  &\Cor\bigg[u_{\operatorname{s}},\frac{\partial u_{\operatorname{s}}}{\partial{\bf n}}\bigg]({\bf x},{\bf x}')
  	= \overline{\Cor\bigg[\frac{\partial u_{\operatorname{s}}}{\partial{\bf n}},u_{\operatorname{s}}\bigg]}({\bf x}',{\bf x})\\
  &\qquad= \E\bigg[\overline{\bigg(\int_{\Gamma({\bf y})}\Phi({\bf x}, {\bf z})
	\frac{\partial u_{\operatorname{s}}}{\partial{\bf n}}({\bf z},\omega)\d\sigma_{\bf z}\bigg)}
	\bigg(\int_{\Gamma({\bf y})}\frac{\partial\Phi({\bf x}', {\bf z})}{\partial {\bf n}_{\bf z}}
		u_{\operatorname{s}}({\bf z},{\bf y})\d\sigma_{{\bf z}}\bigg)\bigg],
\end{align*}
we have for any ${\bf x},{\bf x}'\not\in B_R({\bf 0})$ the deterministic expression 
\begin{equation}\label{eq:cor}
\begin{aligned}
  \Cor[u_{\operatorname{s}}]({\bf x},{\bf x}') = \int_\Sigma\int_\Sigma
  	\bigg\{&\Phi({\bf x}, {\bf z})\overline{\Phi({\bf x}', {\bf z}')}\Cor\!\!\bigg[\frac{\partial u_{\operatorname{s}}}{\partial{\bf n}}\bigg]({\bf z},{\bf z}')\\
  &\ + \Phi({\bf x}, {\bf z})\overline{\frac{\partial\Phi({\bf x}', {\bf z}')}{\partial {\bf n}_{{\bf z}'}}}
  	\Cor\!\!\bigg[\frac{\partial u_{\operatorname{s}}}{\partial{\bf n}},u_{\operatorname{s}}\bigg]({\bf z},{\bf z}')\\
  &\ + \frac{\partial\Phi({\bf x}, {\bf z})}{\partial {\bf n}_{\bf z}}\overline{\Phi({\bf x}', {\bf z}')}
  	\Cor\!\!\bigg[u_{\operatorname{s}},\frac{\partial u_{\operatorname{s}}}{\partial{\bf n}}\bigg]({\bf z},{\bf z}')\\
  &\ + \frac{\partial\Phi({\bf x}, {\bf z})}{\partial {\bf n}_{\bf z}}
	\overline{\frac{\partial\Phi({\bf x}', {\bf z}')}{\partial {\bf n}_{{\bf z}'}}}
	\Cor[u_{\operatorname{s}}]({\bf z},{\bf z}')\bigg\}\d\sigma_{{\bf z}'}\d\sigma_{\bf z}.
\end{aligned}
\end{equation}
As we will see in the next section, this expression can efficiently be computed
if a low-rank approximation of the Cauchy datas' correlations is available.
\section{Numerical results}\label{sec:results}
\subsection{Random scatterer}
For our numerical experiments, we shall
consider a kite-shaped scatterer as nominal 
obstacle, described by the parametrization
\begin{equation}\label{kiteBoundary}
	\boldsymbol\gamma\colon[0,2\pi]\to\Gamma\subset\mathbb{R}^2,\quad
		\phi\mapsto\boldsymbol\gamma(\phi)\isdef
			\begin{bmatrix}
				5\cos(\phi)-3.25\cos(2\phi) \\
				7.5\sin(\phi)\end{bmatrix}.
\end{equation}
The random boundary is then defined in accordance with
\begin{equation}\label{randomPara}
\boldsymbol\gamma(\phi,{\bf y}) = \overline{\boldsymbol\gamma}(\phi)
	+ r(\phi,{\bf y})\begin{bmatrix}\cos(\phi)\\\sin(\phi)\end{bmatrix},
\end{equation}
where $\overline{\boldsymbol\gamma}(\phi)$ denotes the 
kite-shaped boundary \eqref{kiteBoundary} and $r(\phi,{\bf y})$ 
is given by the Fourier series
\begin{equation}\label{randomRadii}
  r(\phi,{\bf y}) = \sum_{k=1}^{\infty}
  	\frac{1}{k^3}\big\{\sin(k\phi) y_{2k-1} + \cos(k\phi) y_{2k}\big\}.
\end{equation}
For the numerical simulation, we truncate this series after
$1000$ terms. 

\begin{figure}[hbt]
\begin{center}
\includegraphics[trim=250 80 190 60,clip,width=0.5\textwidth]{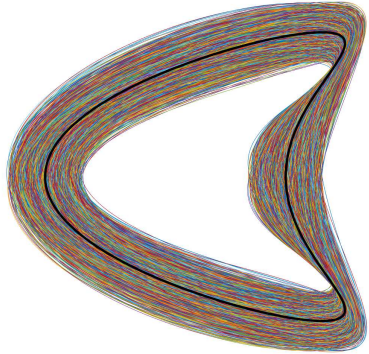}
\caption{\label{fig:RandomDomain}The kite-shaped boundary
(thick black line) and 1000 random perturbations (in colour).}
\end{center}
\end{figure}

Notice that the decay of the coefficients of the random 
fluctuations \eqref{randomRadii} are at the limit case. It 
would hold $r(\cdot,{\bf y})\in C_{\per}^2([0,2\pi])$ if the decay 
of the series $\{a_k\}_k$ was just a bit stronger. A visualization 
of $1000$ samples of this boundary is found in 
Figure~\ref{fig:RandomDomain}.

\subsection{Statistics at the artificial interface}
For the numerical solution of the boundary integral equation 
\eqref{eq:BIE1}, we apply the Nystr\"om method to discretize 
the acoustic single and double layer operators. 
Given the 
parametrization \eqref{randomPara} for a specific instance 
${\bf y}\in[-1,1]^{1000}$, the 
method applies the trapezoidal rule in the $n = 1000$ equidistantly 
distributed points $\varphi_i = \nicefrac{2\pi i}{n}$, $i=1,\ldots,n$, and 
is along the lines of \cite[Chapter 12]{K}. An appropriate desingularization 
technique based on trigonometric Lagrange polynomials is 
employed to deal with the singularities of the acoustic single 
and double layer operators. We remark that this method converges exponentially 
provided that the boundary under consideration is analytical.
We refer the reader to \cite[Chapter 12]{K} for all the details. 

We likewise subdivide the artificial interface 
\(\Sigma = \partial B_R({\bf 0})\)
by $n = 1000$ equidistantly distributed points 
\[
  {\bf z}_j = \big[R\cos(\nicefrac{2\pi j}{n}),R\sin(\nicefrac{2\pi j}{n})\big]^\top,\quad j = 1,\ldots,n. 
\]
In these points, we compute the expectations $\E[u_{\operatorname{s}}]({\bf z}_j)$
and $\E[\nicefrac{\partial u_{\operatorname{s}}}{\partial {\bf n}}]({\bf z}_j)$ in accordance 
with \eqref{eq:expectation3} and \eqref{eq:expectation4}, respectively.
To that end, we employ the quasi-Monte Carlo method based on 
$10\,000$ Halton points, cf.~\cite{NIE}. In \cite{Wan02} it is shown that this
quadrature converges independent of the parameter dimension, if the derivatives with
respect to the parameter \({\bf y}\in\square\) decay sufficiently fast. \red{A thorough analysis
of the computational work of the Nystr{\"o}m method in combination with the quasi-Monte Carlo
quadrature can be found in \cite{GP18}.}

In addition to the expectations of the Cauchy data, we compute the 
corresponding two-point correlation matrix
\begin{equation}\label{eq:C}
  {\bf C} = \begin{bmatrix} {\bf C}_{1,1} & {\bf C}_{1,2} \\
  	{\bf C}_{1,2}^\star & {\bf C}_{2,2}\end{bmatrix}\in\mathbb{C}^{2n\times 2n},
\end{equation}
where
\[
  {\bf C}_{1,1} := \big[\Cor[u_{\operatorname{s}}]({\bf z}_j,{\bf z}_{j'})\big]_{j,j'=1}^n,\quad
  {\bf C}_{2,2} := \Bigg[\Cor\bigg[\frac{\partial u_{\operatorname{s}}}{\partial {\bf n}}\bigg]({\bf z}_j,{\bf z}_{j'})\bigg]_{j,j'=1}^n
\]
and
\[
  {\bf C}_{1,2} := \Bigg[\Cor\bigg[u_{\operatorname{s}},\frac{\partial u_{\operatorname{s}}}{\partial {\bf n}}\bigg]({\bf z}_j,{\bf z}_{j'})\Bigg]_{j,j'=1}^n.
\]

\subsection{Low-rank approximation of the two-point correlation}
While the computation of $\E[u_{\operatorname{s}}]({\bf x})$ at a point ${\bf x}\not\in B_R({\bf 0})$ 
by \eqref{eq:expectation2} is straightforward,
the computation of the variance $\V[u_{\operatorname{s}}]({\bf x})$ in accordance 
with \eqref{eq:var} amounts to the computation of $\Cor[u_{\operatorname{s}}]
({\bf x},{\bf x})$. This requires the approximation of the double 
integral over $\Sigma$. We apply again the trapezoidal rule, 
having thus to evaluate
\begin{align*}
  \Cor[u_{\operatorname{s}}]({\bf x},{\bf x})\approx \frac{1}{(2R\pi n)^2}\sum_{j,j'=1}^n\bigg\{
  &\Phi({\bf x},{\bf z}_j)\overline{\Phi({\bf x},{\bf z}_{j'})}\Cor[u_{\operatorname{s}}]({\bf z}_j,{\bf z}_{j'})\\
 &\ + \Phi({\bf x},{\bf z}_j)\overline{\frac{\partial\Phi({\bf x},{\bf z}_{j'})}{\partial {\bf n}_{{\bf z}_{j'}}}}
    \Cor\bigg[u_{\operatorname{s}},\frac{\partial u_{\operatorname{s}}}{\partial {\bf n}}\bigg]({\bf z}_j,{\bf z}_{j'})\\
     &\ + \overline{\Phi({\bf x},{\bf z}_j)}\frac{\partial\Phi({\bf x},{\bf z}_{j'})}{\partial {\bf n}_{{\bf z}_{j'}}}
    \overline{\Cor\bigg[u_{\operatorname{s}},\frac{\partial u_{\operatorname{s}}}{\partial {\bf n}}\bigg]({\bf z}_{j'},{\bf z}_j)}\\
 &\ + \frac{\partial\Phi({\bf x},{\bf z}_j)}{\partial {\bf n}_{{\bf z}_j}}
   \overline{\frac{\partial\Phi({\bf x},{\bf z}_{j'})}{\partial {\bf n}_{{\bf z}_{j'}}}}
    \Cor\bigg[\frac{\partial u_{\operatorname{s}}}{\partial {\bf n}}\bigg]({\bf z}_j,{\bf z}_{j'})\bigg\}.
\end{align*}
The respective evaluations of the two-point correlation 
functions of the Cauchy data at $\Sigma$ are stored in 
the matrix ${\bf C}$ from \eqref{eq:C}. We conclude that
the cost of a naive evaluation scales quadratically 
in the number of degrees of freedom used at the artificial 
interface $\Sigma$.

In order to speed-up the computations if the variance $\V[u_{\operatorname{s}}]({\bf x})$
has to be computed in many points, we propose to compute first a
low-rank approximation of the two-point correlation function of the
Cauchy data at $\Sigma$. In accordance with \cite{HPS}, we 
apply the pivoted Cholesky decomposition to get a low-rank 
approximation
\begin{equation}\label{eq:PCD}
  {\bf C} \approx {\bf LL}^\star = \sum_{i=1}^m \boldsymbol\ell_i\boldsymbol\ell_i^\star
\end{equation}
where ${\bf L}=[\boldsymbol\ell_1,\ldots,\boldsymbol\ell_m]\in
\mathbb{C}^{2n\times m}$ with $m\le n$. \red{Note that the truncation 
error can rigorously be controlled in terms of the trace. Hence, the
pivoted Cholesky decomposition is truncated if
\begin{equation}\label{eq:PCtrunc}
\operatorname{trace}({\bf C}-{\bf LL}^\star)
<\varepsilon\operatorname{trace}({\bf C})
\end{equation}
for some \(\varepsilon>0\). For all the details, we refer to \cite{HPS,HPS15}.
We remark that we would still end up with a separable expansion if
the covariance of the Cauchy data did not admit a low-rank representation.}

Having the low-rank approximation \eqref{eq:PCD} at hand, 
we arrive at
\begin{equation}\label{eq:low-rank-cor}
\begin{aligned}
  \Cor[u_{\operatorname{s}}]({\bf x},{\bf x}) &\approx \frac{1}{(2R\pi n)^2}\sum_{i=1}^m
  	\Bigg|\sum_{j=1}^n\bigg[\Phi({\bf x},{\bf z}_j)\ell_{i,j}
		+\frac{\partial\Phi({\bf x},{\bf z}_j)}{\partial {\bf n}_{{\bf z}_j}}\ell_{i,n+j}\bigg]\Bigg|^2.
\end{aligned}
\end{equation}
Therefore, the evaluation of $ \Cor[u_{\operatorname{s}}]({\bf x},{\bf x}')$
requires only $\mathcal{O}(nm)$ operations instead of
$\mathcal{O}(n^2)$ operations. If $m\ll n$, this reduces
the computational cost considerably, especially
since $m$ depends only on the desired accuracy and
thus only weakly on $n$. 

\begin{table}[hbt]
\begin{center}
\begin{tabular}{|c|ccccc|}\hline
 \multicolumn{6}{|c|}{rank of the low-rank approximation}\\\hline
 $R$ & $\kappa = 1$ & $\kappa = 2$ & $\kappa = 4$ & $\kappa = 8$ & $\kappa = 16$ \\\hline
 11 & 48 & 56 & 85 & 131 & 193 \\
 12 & 39 & 51 & 83 & 131 & 194 \\
 13 & 35 & 49 & 84 & 132 & 195 \\ 
 14 & 32 & 49 & 83 & 131 & 195 \\
 15 & 31 & 49 & 84 & 132 & 194 \\ \hline
\end{tabular}
\caption{\label{tab:m}Ranks $m$ of the low-rank approximation
of the two-point correlation of the Cauchy data at $\Sigma$ for
varying radius $R$ and wavenumber $\kappa$.}
\end{center}
\end{table}

\begin{figure}[htb]
\begin{center}
\begin{minipage}{0.48\textwidth}
\begin{center}
\(\kappa=1\)
\includegraphics[trim = 0 0 0 0, clip,width=\textwidth]{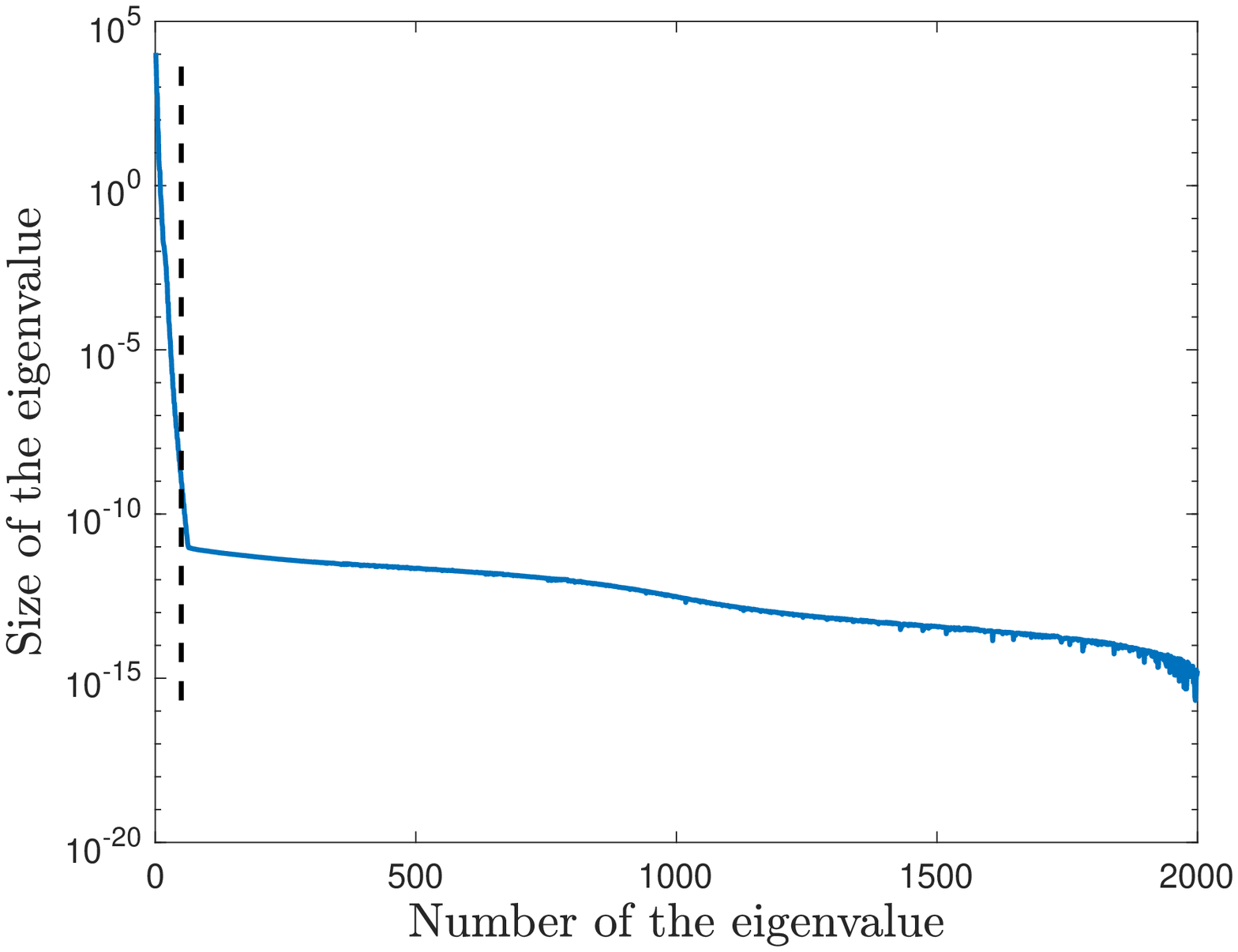}
\end{center}
\end{minipage}\hfill
\begin{minipage}{0.48\textwidth}
\begin{center}
\(\kappa=2\)
\includegraphics[trim = 0 0 0 0, clip,width=\textwidth]{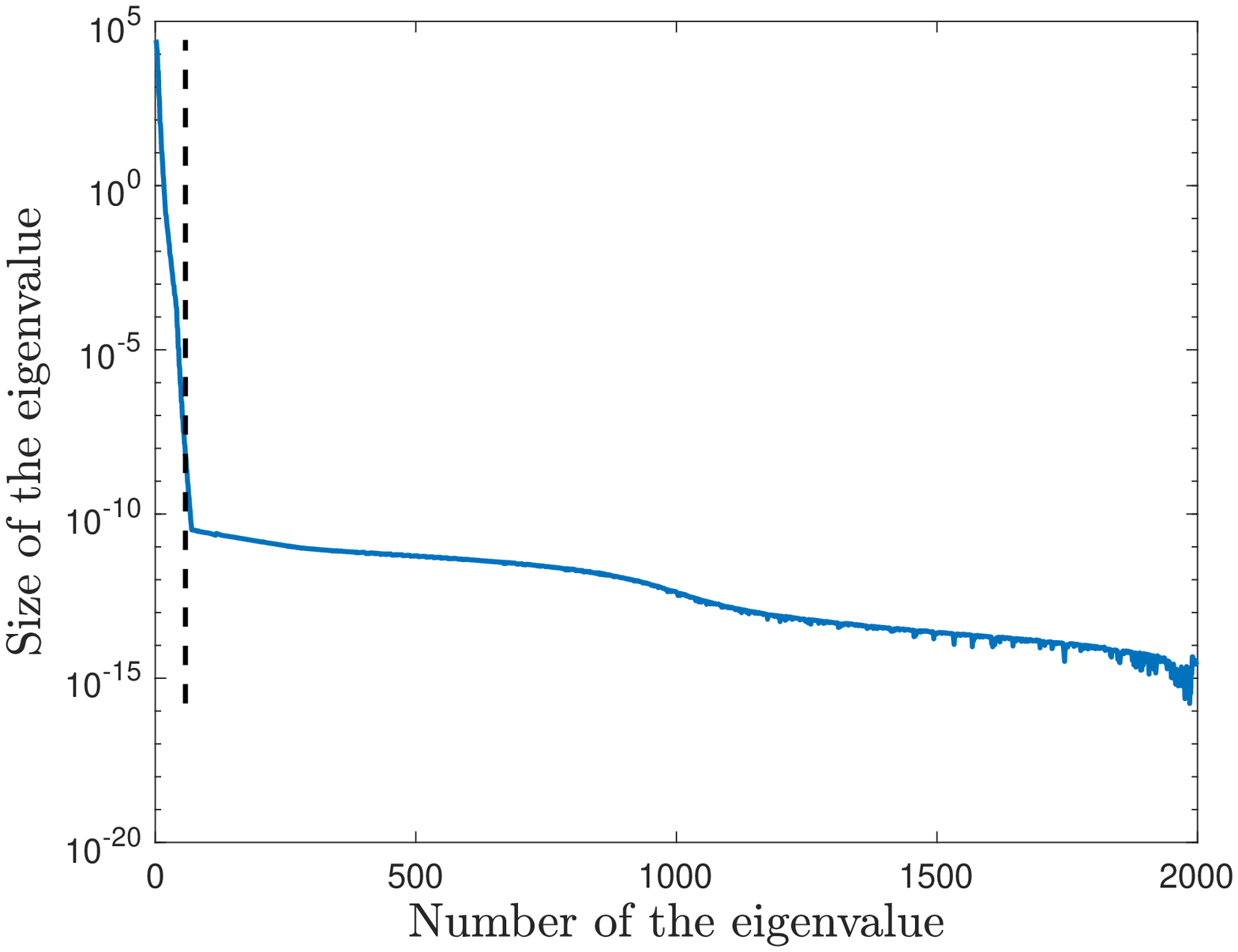}
\end{center}
\end{minipage}\\[1em]
\begin{minipage}{0.48\textwidth}
\begin{center}
\(\kappa=4\)
\includegraphics[trim = 0 0 0 0, clip,width=\textwidth]{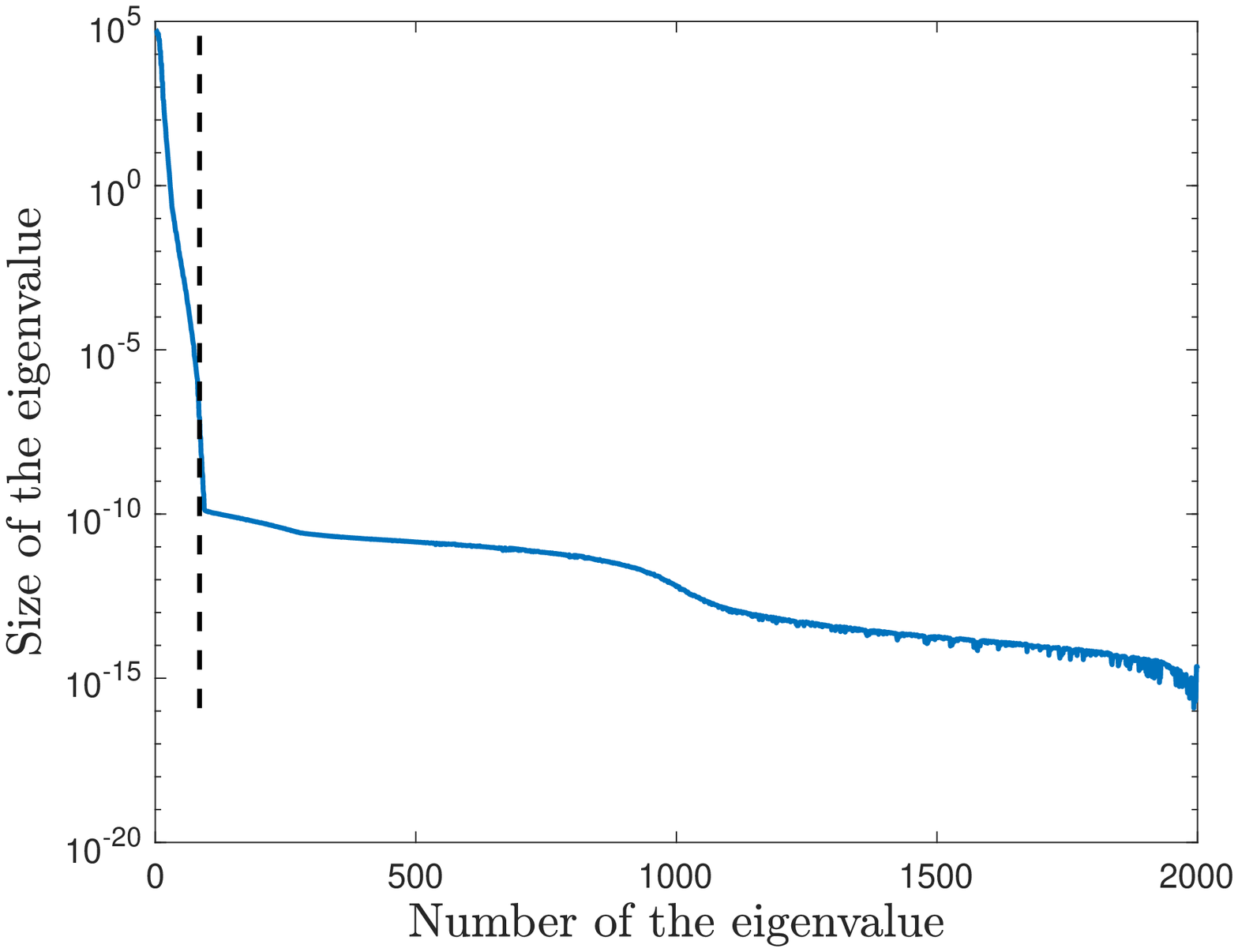}
\end{center}
\end{minipage}\hfill
\begin{minipage}{0.48\textwidth}
\begin{center}
\(\kappa=8\)
\includegraphics[trim = 0 0 0 0, clip,width=\textwidth]{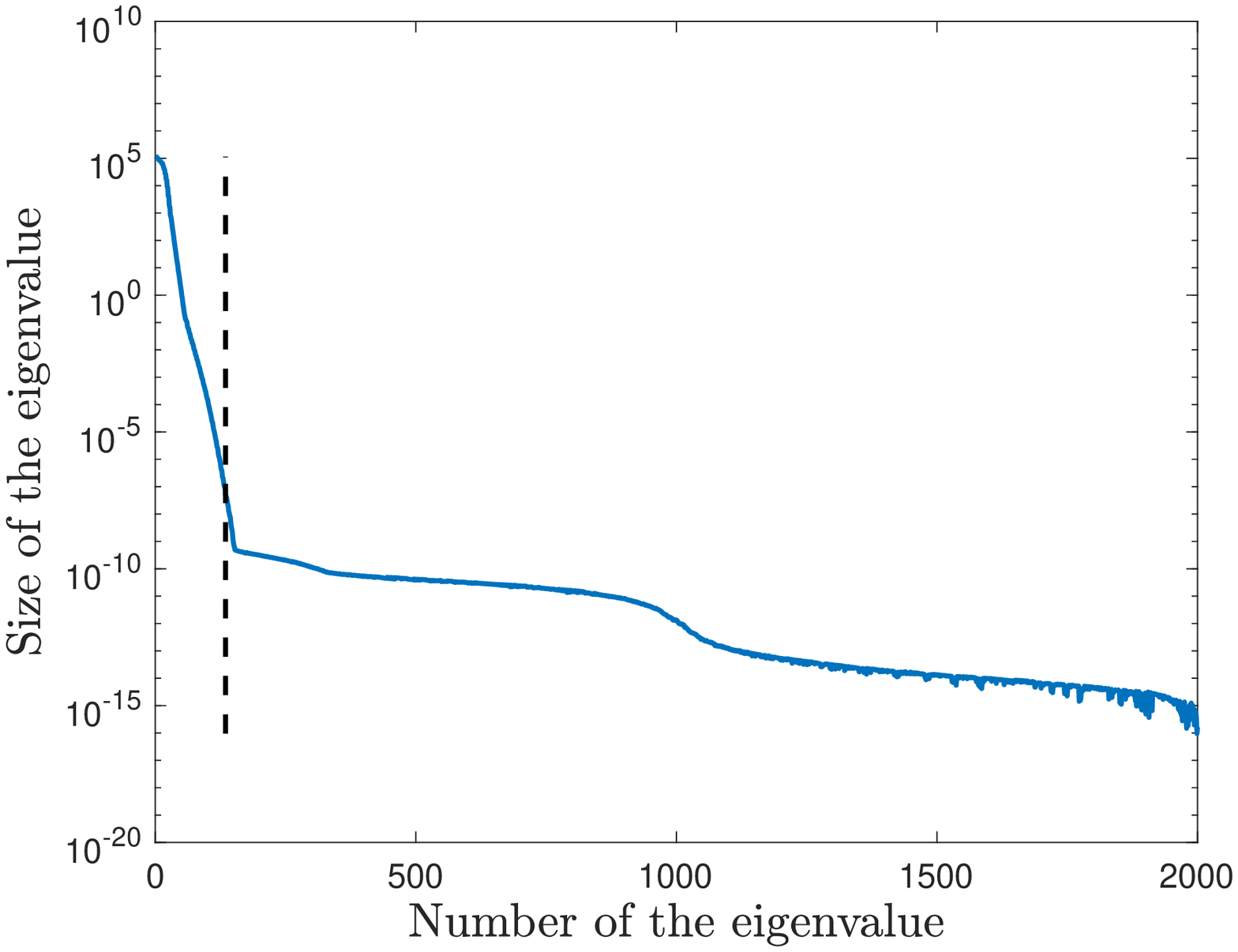}
\end{center}
\end{minipage}
\caption{\label{fig:decayEV} Decay of the eigenvalues of the 
covariance matrix \({\bf C}\) for the wavenumbers \(\kappa=1,2,4,8\) and \(R=11\).}
\end{center}
\end{figure}

In order to demonstrate the efficiency of the low-rank 
approximation, we consider again the randomly 
perturbed kite-shaped scatterer, given by \eqref{randomPara} 
and \eqref{randomRadii}. The radius of the artificial interface 
is varying in accordance with $R = 11,12,\ldots,15$ and the 
wavenumber is varying in accordance with $\kappa = 1,2,4,8,16$.
The number of equidistant points on $\Sigma$ is 1000
and the number of boundary elements on $\Gamma({\bf y})$
is also 1000. Note that the incident wave has been chosen to come 
from the left, i.e., ${\bf d} = [1,0]^\top$, and the 
\red{upper bound for the relative truncation
error of the Cauchy data's pivoted Cholesky decomposition is $10^{-12}$,
cp.\ \eqref{eq:PCtrunc}}.	
The corresponding results are found in Table~\ref{tab:m}. As can be seen,
the pivoted cholesky decomposition converges very rapidly, where the
determined rank decreases for increasing \(R\). In order to provide a
better intuition of the Cauchy data's covariance, we have also depicted
the corresponding eigenvalues for \(R=11\) in Figure~\ref{fig:decayEV}.

\subsection{Scattered field computation}\label{sec:scatteredfield}
\begin{figure}[htb]
\begin{center}
\begin{tabular}{c|c|c|c}
\(\kappa\) & nominal & expectation & variance\\
\hline & & & \\[-1.1em]
\raisebox{3.4em}{\(\kappa=1\)} &
\includegraphics[trim = 65 25 65 25, clip,width=0.27\textwidth]{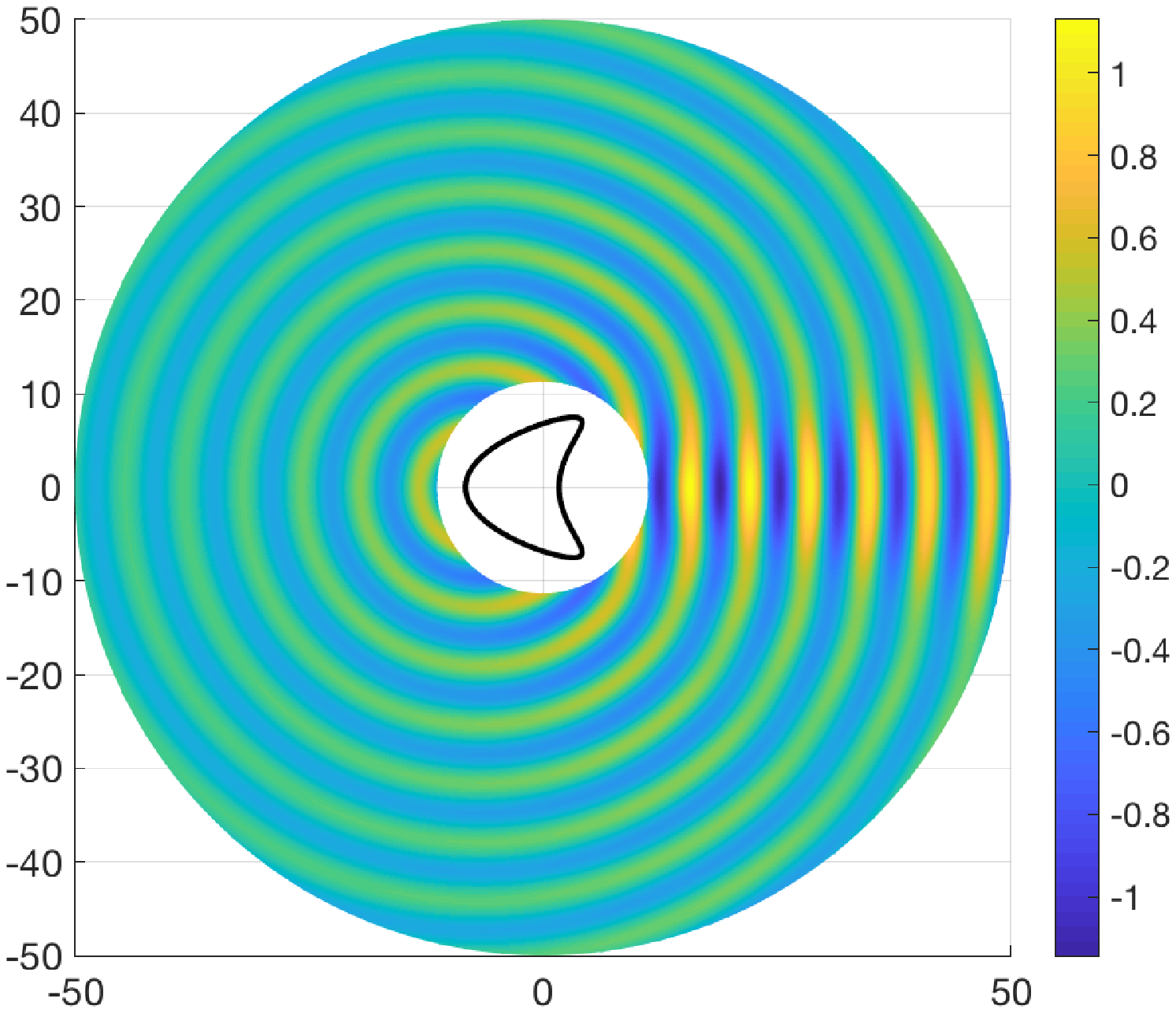} &
\includegraphics[trim = 65 25 65 25, clip,width=0.27\textwidth]{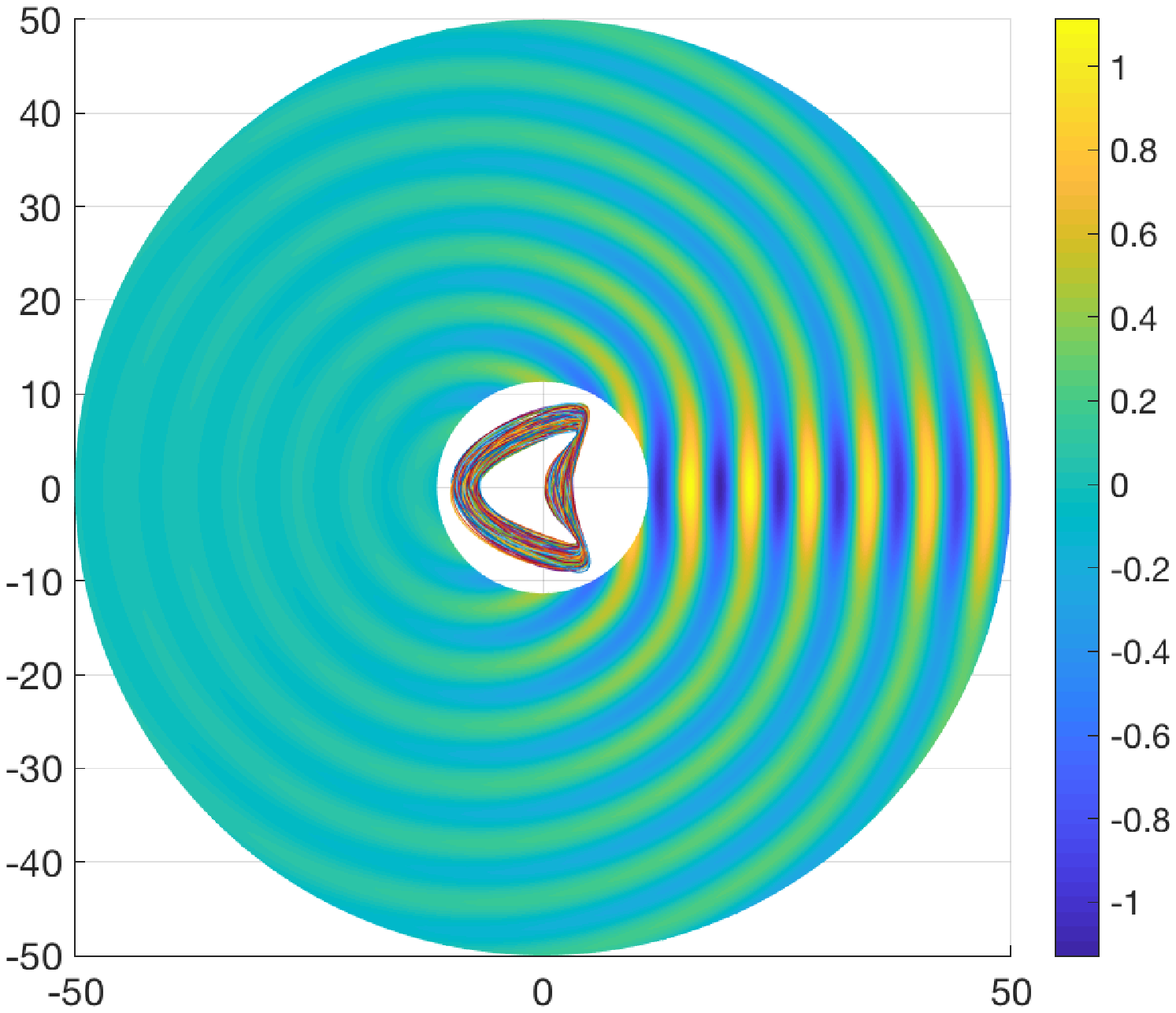}&
\includegraphics[trim = 65 25 65 25, clip,width=0.27\textwidth]{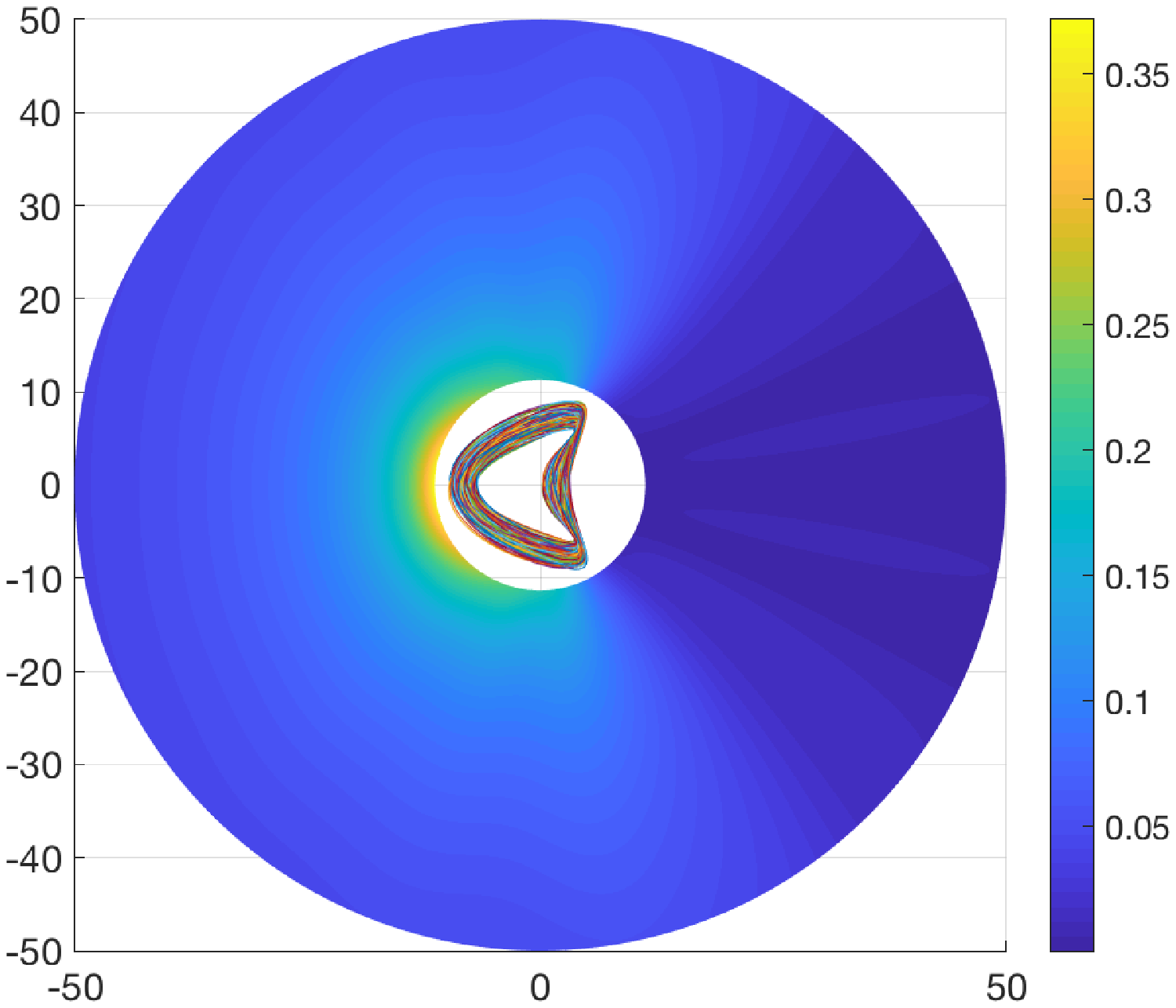}\\[0em]
\hline & & & \\[-1.1em]
\raisebox{3.4em}{\(\kappa=2\)}&
\includegraphics[trim = 65 25 65 25, clip,width=0.27\textwidth]{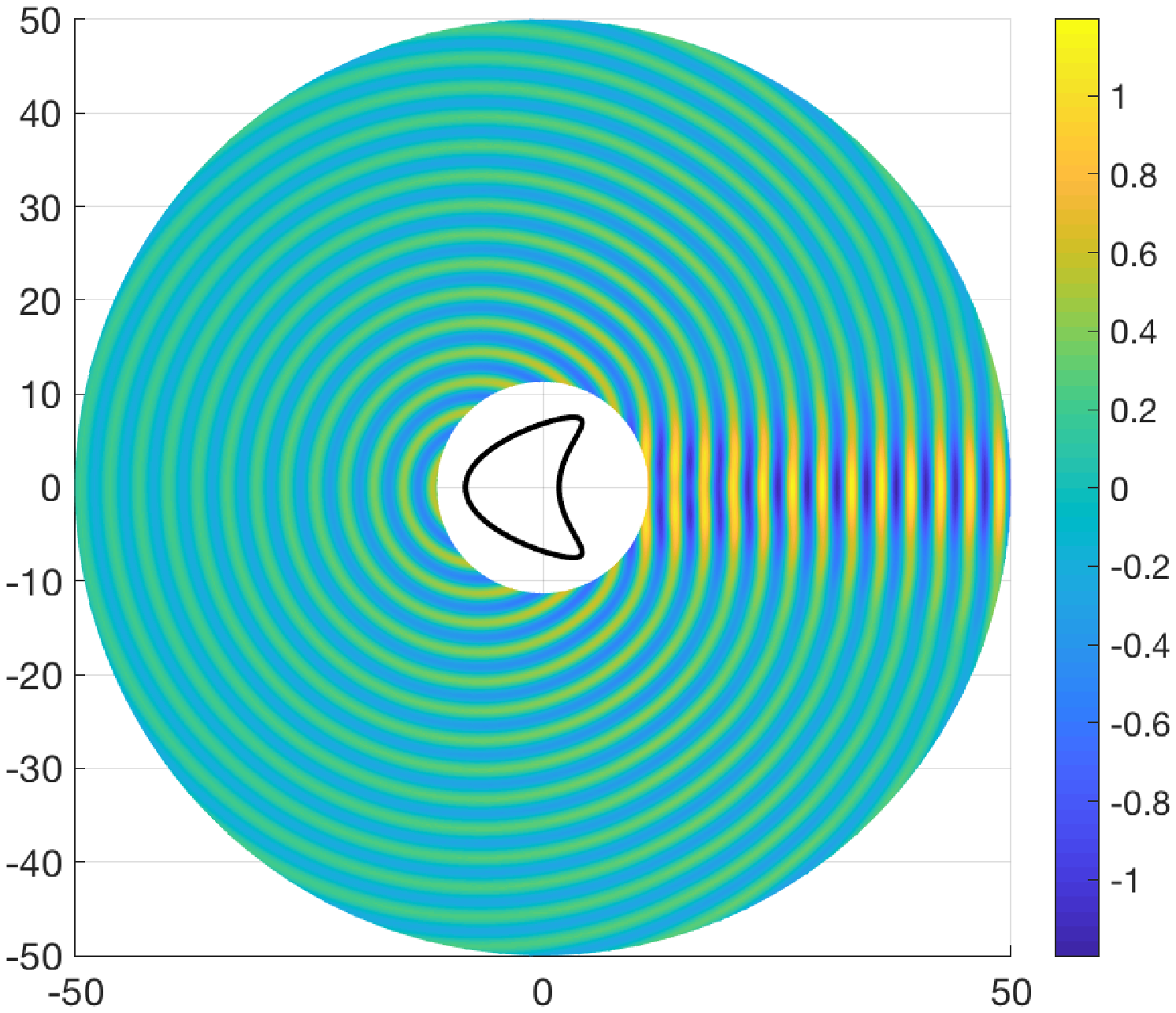} &
\includegraphics[trim = 65 25 65 25, clip,width=0.27\textwidth]{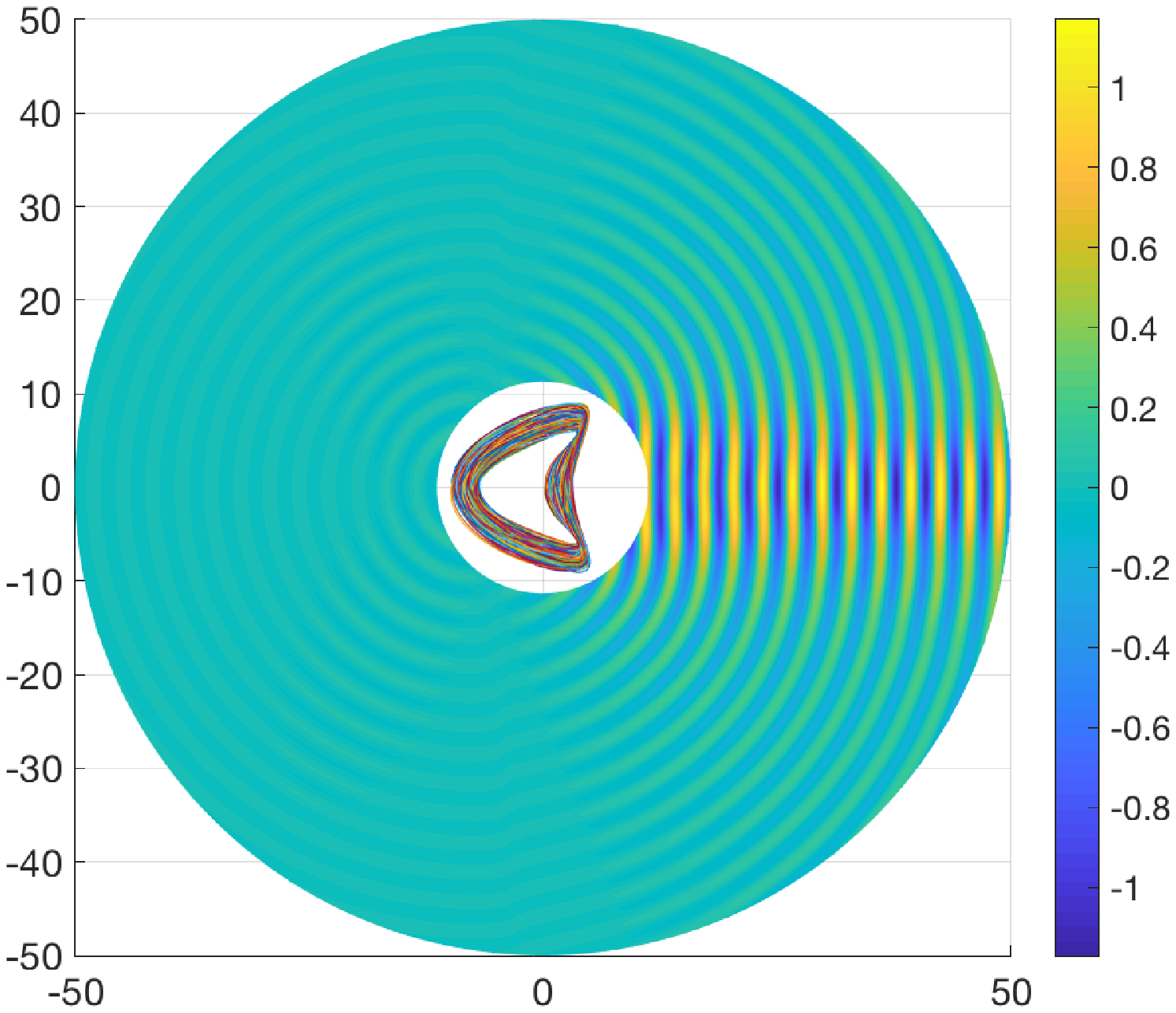} &
\includegraphics[trim = 65 25 65 25, clip,width=0.27\textwidth]{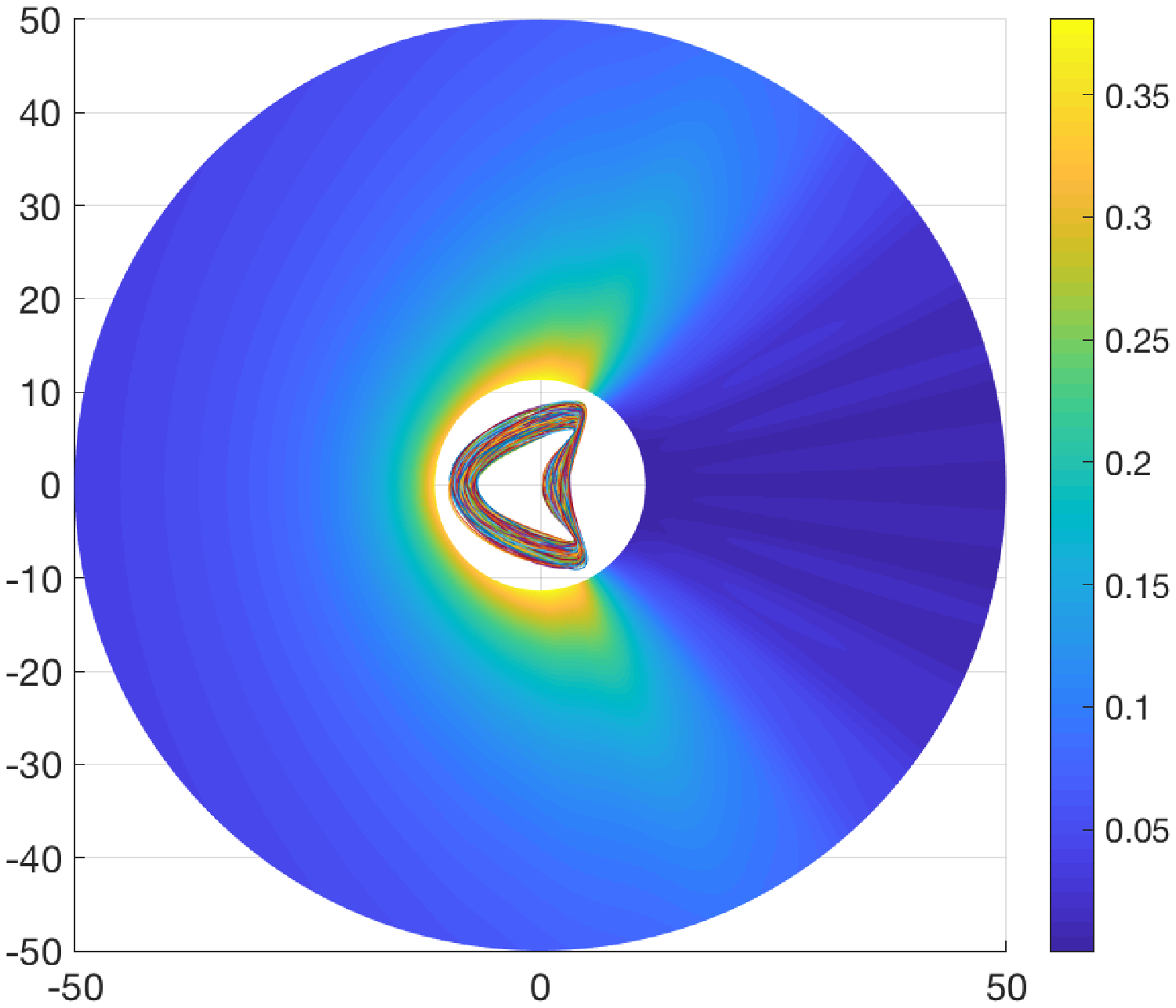}\\[0em] \hline & & & \\[-1.1em]
\raisebox{3.4em}{\(\kappa=4\)}&
\includegraphics[trim = 65 25 65 25, clip,width=0.27\textwidth]{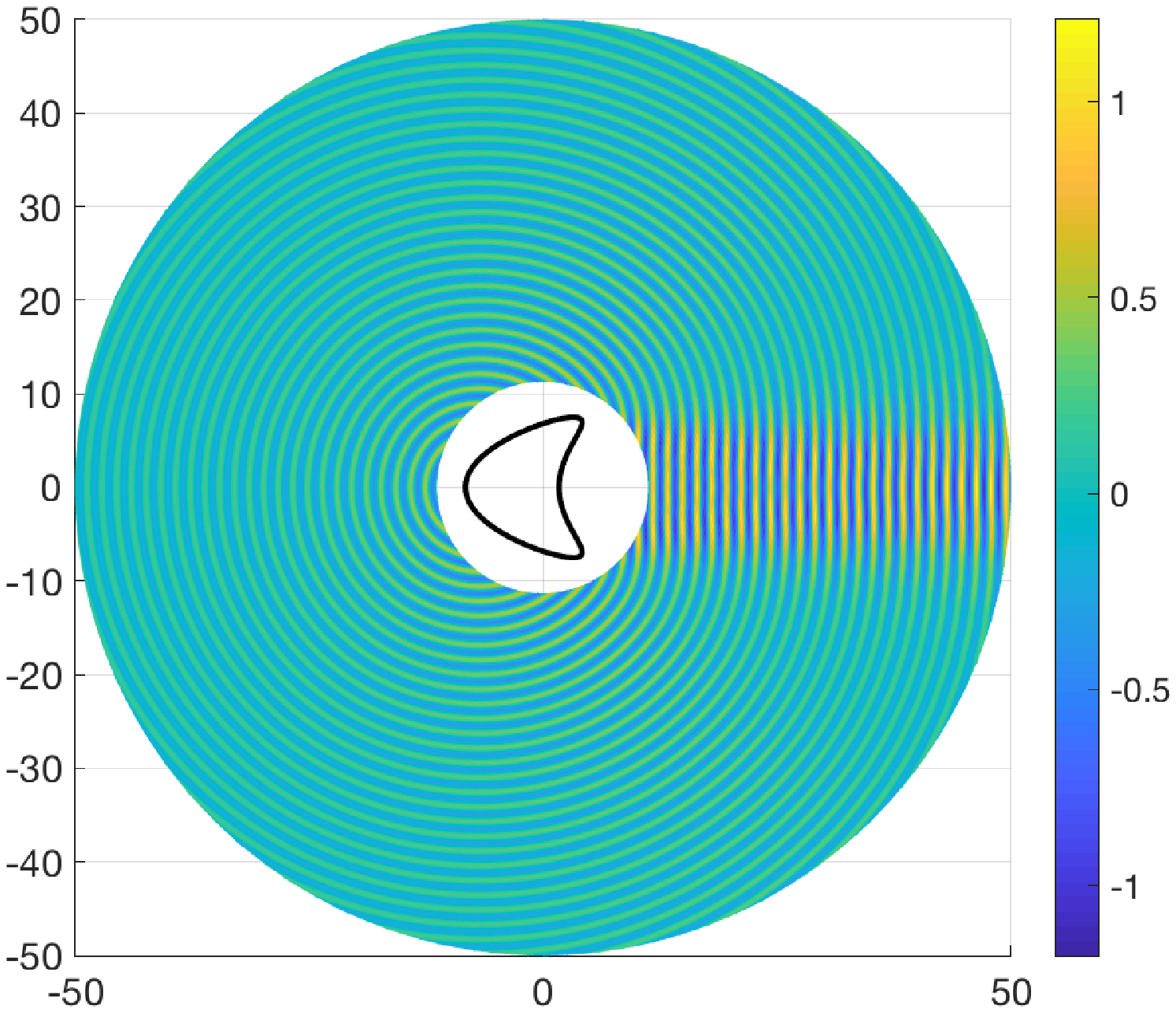} &
\includegraphics[trim = 65 25 65 25, clip,width=0.27\textwidth]{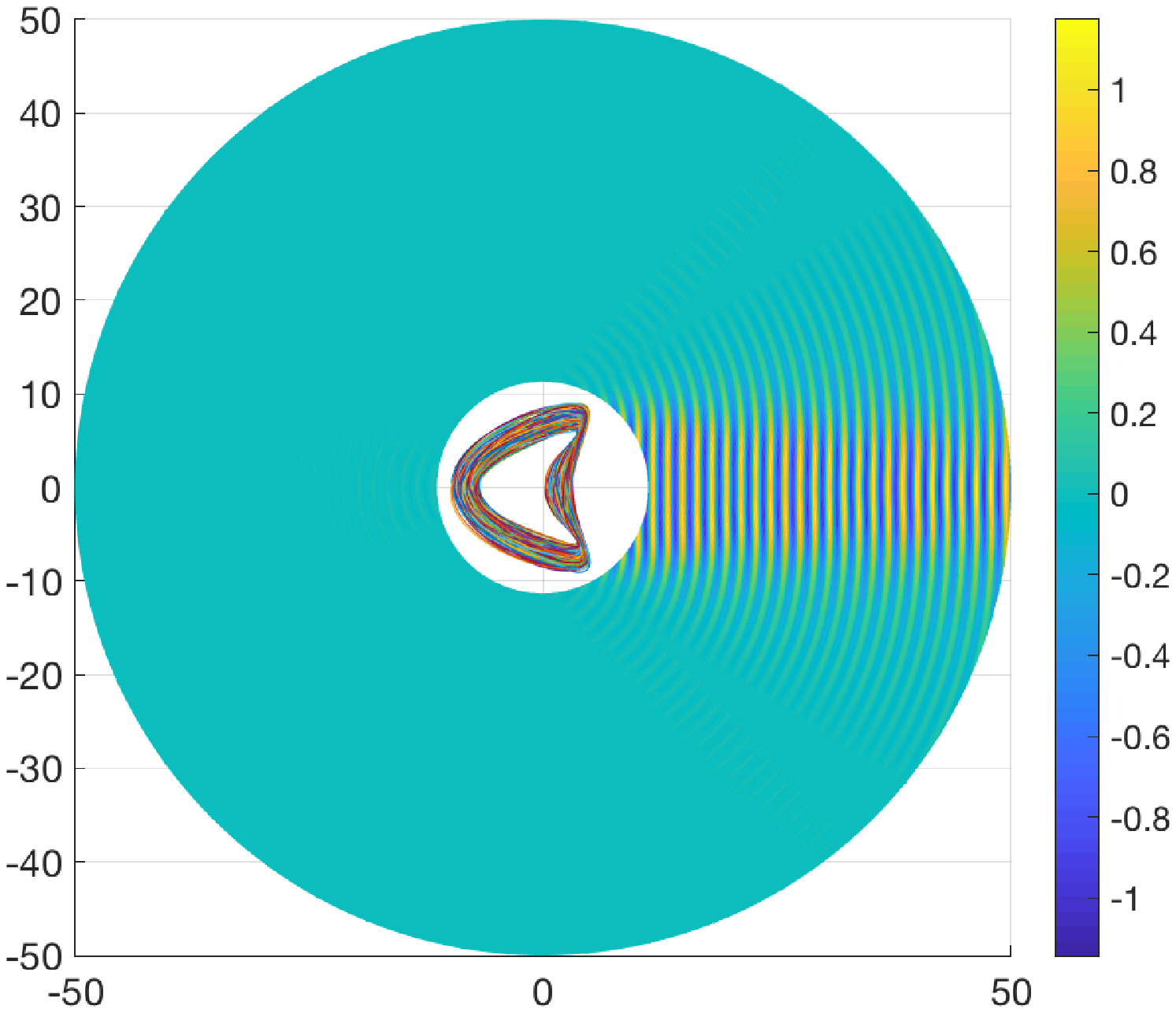} &
\includegraphics[trim = 65 25 65 25, clip,width=0.27\textwidth]{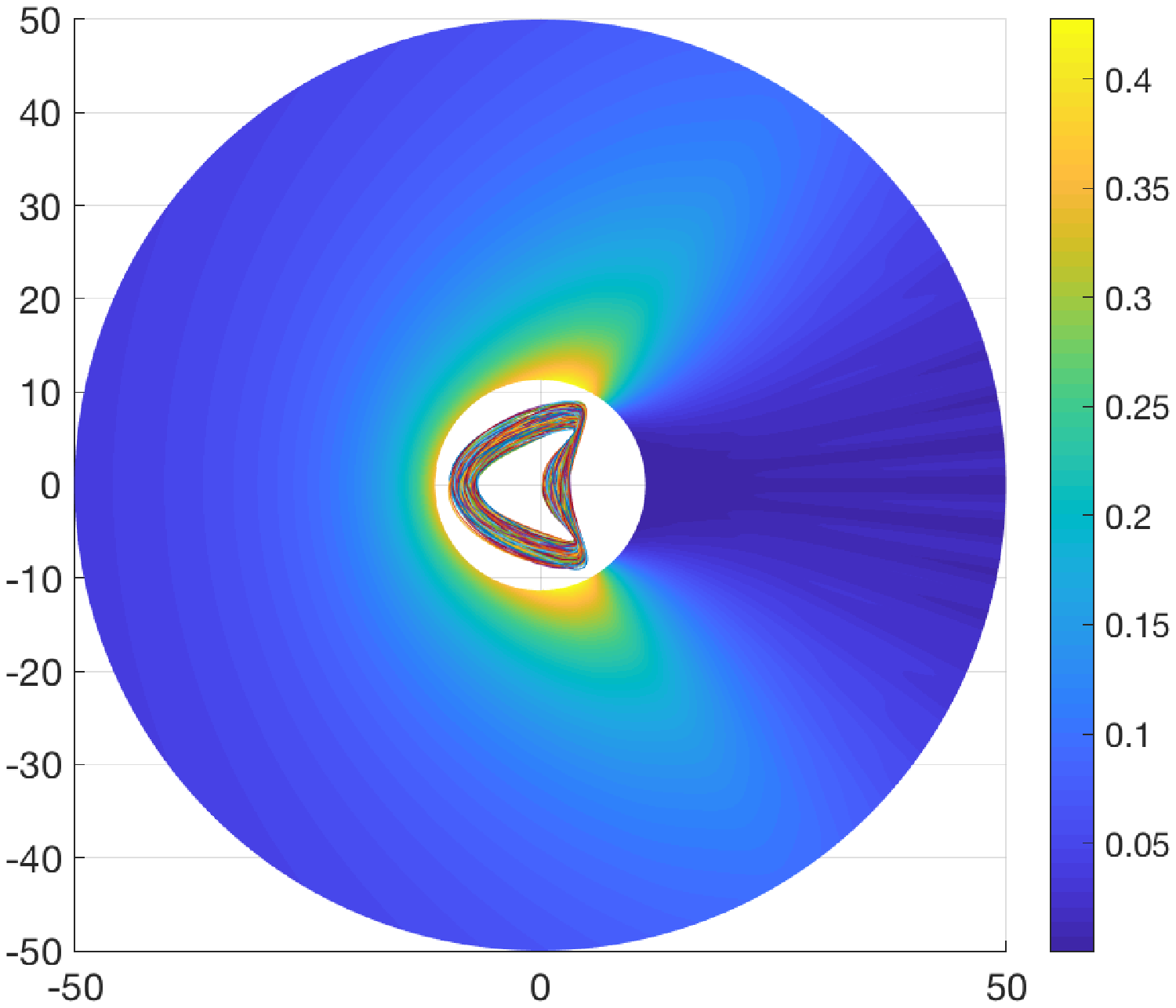}\\[0em] \hline & & & \\[-1.1em]
\raisebox{3.4em}{\(\kappa=8\)}&
\includegraphics[trim = 65 25 65 25, clip,width=0.27\textwidth]{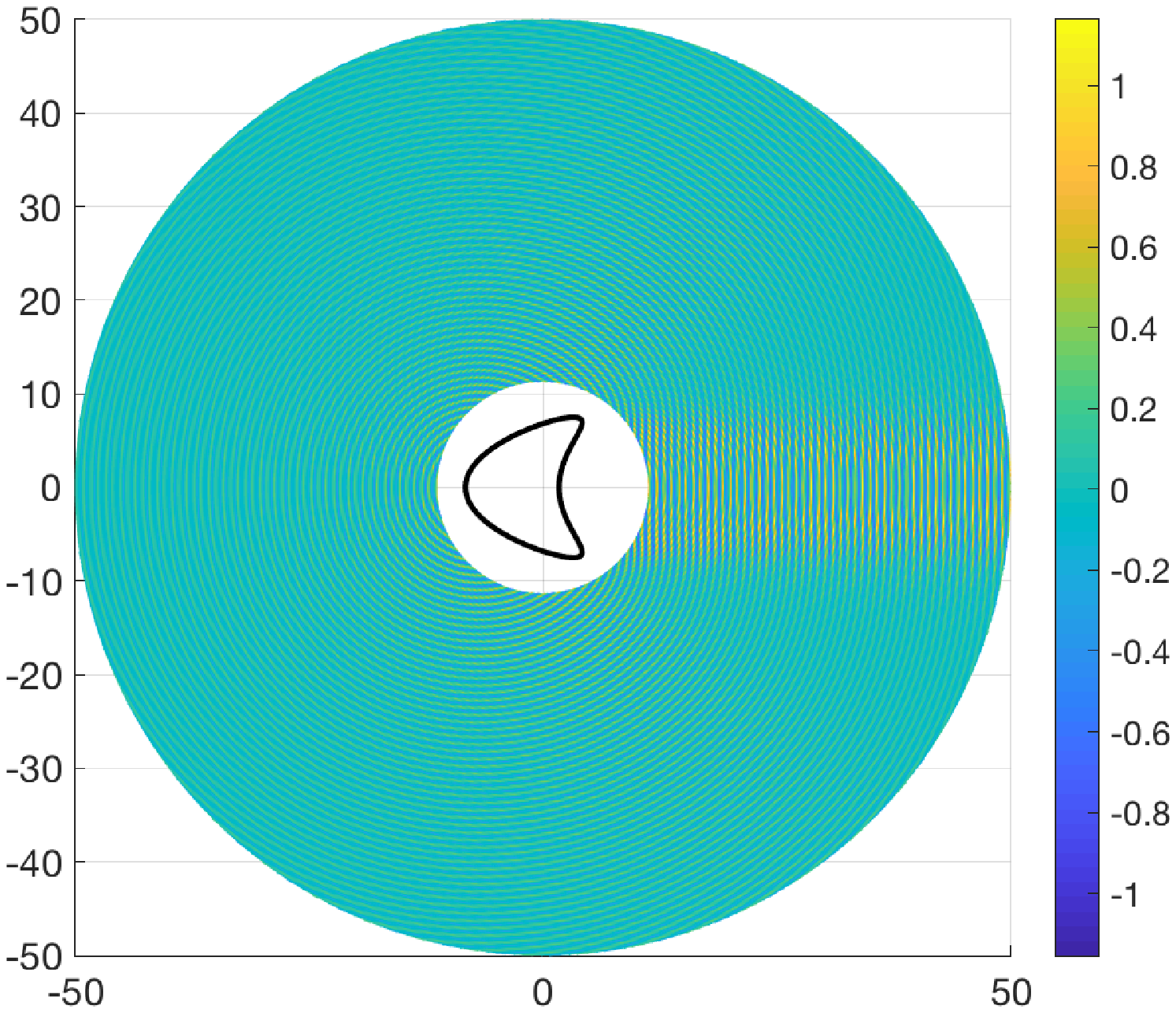} &
\includegraphics[trim = 65 25 65 25, clip,width=0.27\textwidth]{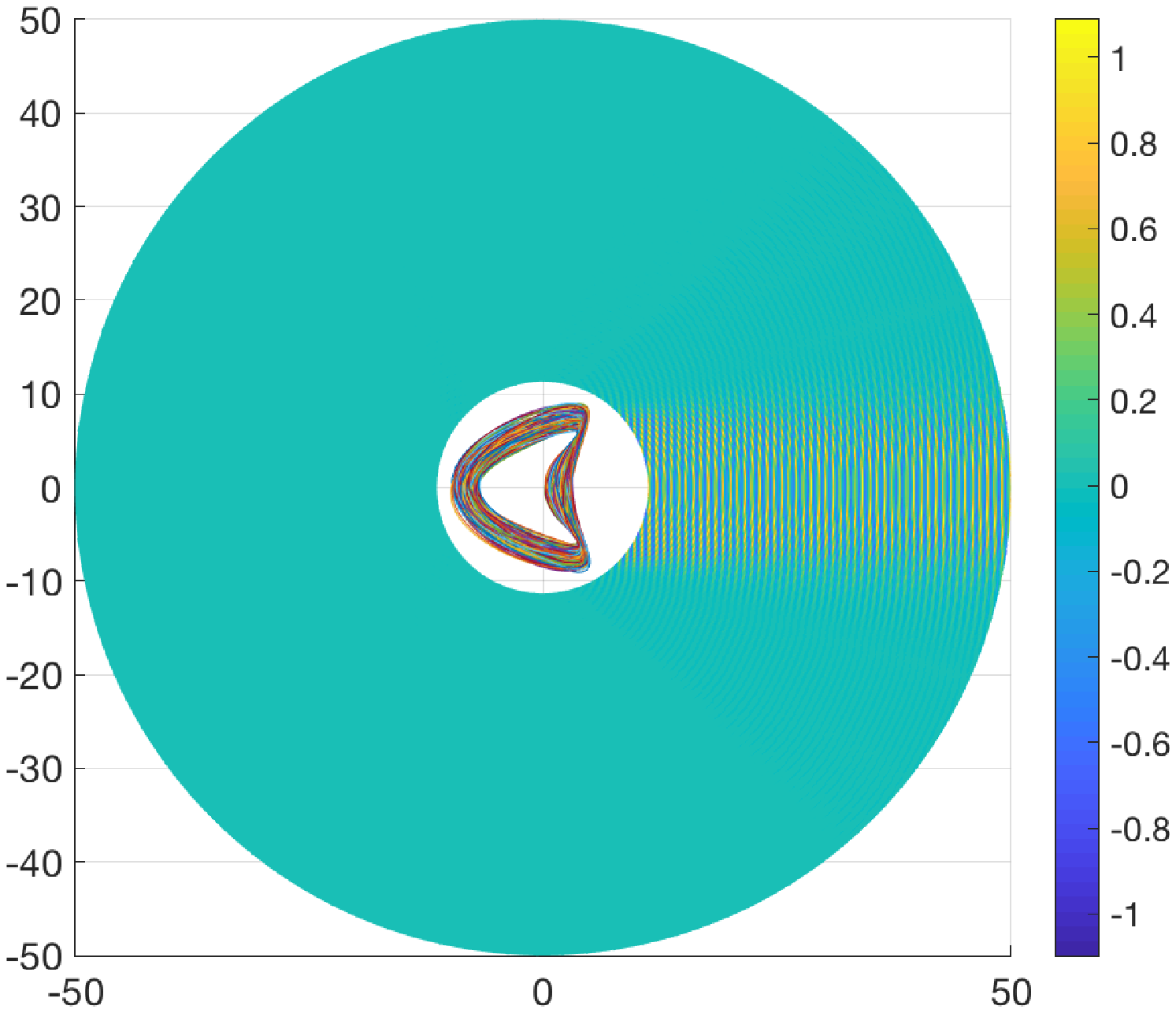} &
\includegraphics[trim = 65 25 65 25, clip,width=0.27\textwidth]{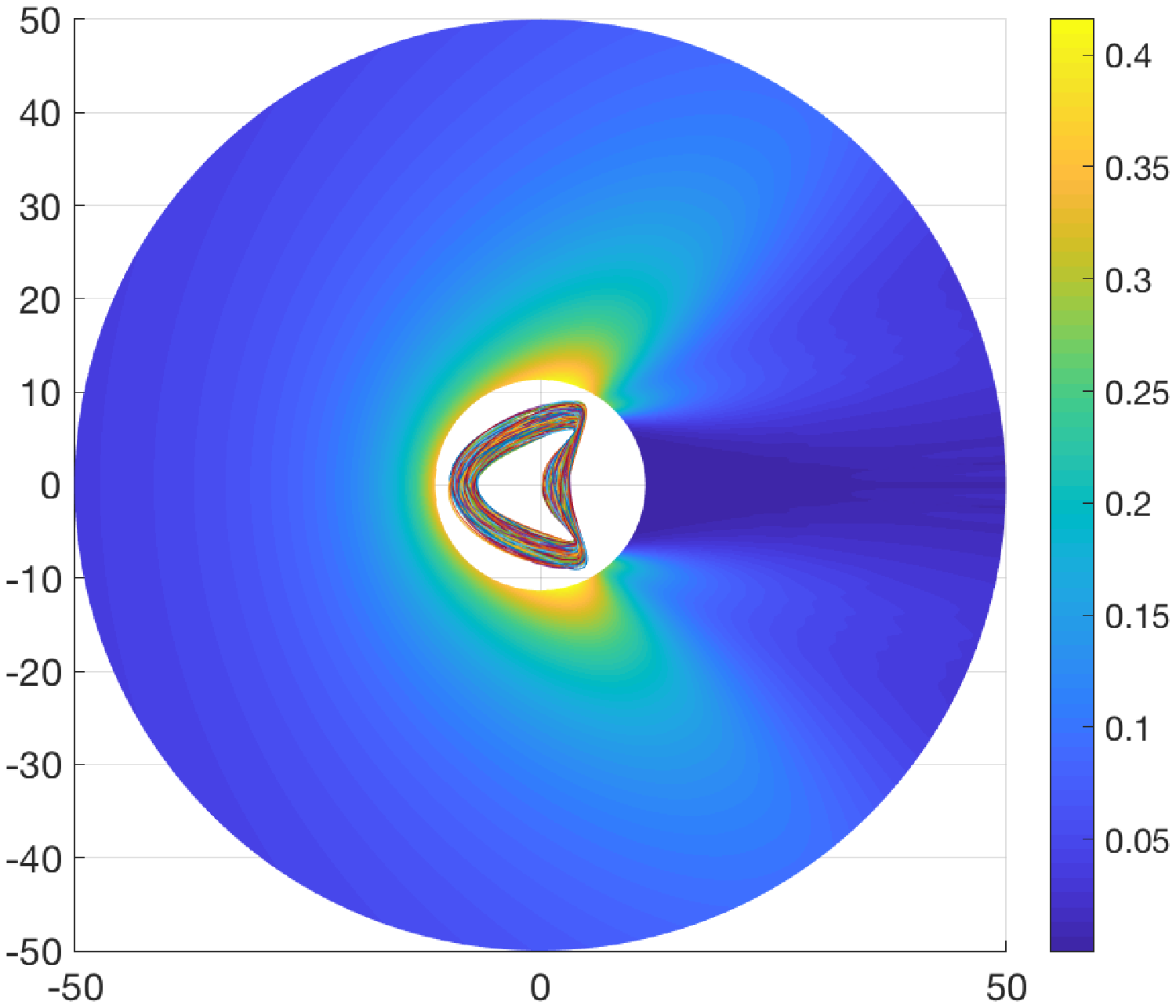}\\[0em]
\end{tabular}
\caption{\label{fig:scattering}The total field of the nominal scatterer (left), 
the expected total field for the random scatterer (middle) and
associated variance (right) for the wavenumbers
$\kappa=1,2,4,8$.} 
\end{center}
\end{figure}

We choose $R=11$ and compute the expectation and variance
of the scattered field on the disc $\{{\bf x}\in\mathbb{R}^2:
R\le\|{\bf x}\|_2\le 50\}$ in accordance with \eqref{eq:expectation2}
and \eqref{eq:var} using \eqref{eq:low-rank-cor}, where the
incident wave comes again from the left, i.e., ${\bf d} = [1,0]^\top$.
The results are found in Figure~\ref{fig:scattering}. For 
comparison, the scattered wave in case of the unperturbed 
kite-shaped scatterer is found in the first column. In the 
second column, the expected total wave is found. 
Finally, the variance of the total wave is found in 
the third column. The rows correspond to the wavenumber: 
the first row corresponds to $\kappa = 1$, the second row
corresponds to $\kappa = 2$, the third row corresponds to 
$\kappa = 4$, and the fourth row corresponds to $\kappa = 8$.

One observes that, compared to the total wave 
of the unperturbed scatterer, the expected total
wave is blurred towards the left, i.e., in directions opposed to the 
direction of the incoming wave. This is caused by the 
different reflections at the perturbed scatterer which 
interfere. In the shadow region, i.e., towards the
right, the expected total wave and the total
wave of the unperturbed scatterer basically coincide.
This observation is also underpinned by the variance 
of the total wave, which is maximal on the left
of the scatterer and nearly 0 in the shadow region. 
Notice that the described smoothing effect becomes 
stronger as the wavenumber increases. 

\subsection{Direction of the incident wave}
\begin{figure}[htb]
\begin{center}
\begin{tabular}{c|c|c|c}
\({\bf d}\) & nominal & expectation & variance\\
\hline & & & \\[-1.1em]
\raisebox{3.4em}{\hspace*{-1em}\({\bf d}=\begin{bmatrix}-1\\ \phantom{-}0\end{bmatrix}\)\hspace*{-0.5em}} &
\includegraphics[trim = 65 25 65 25, clip,width=0.27\textwidth]{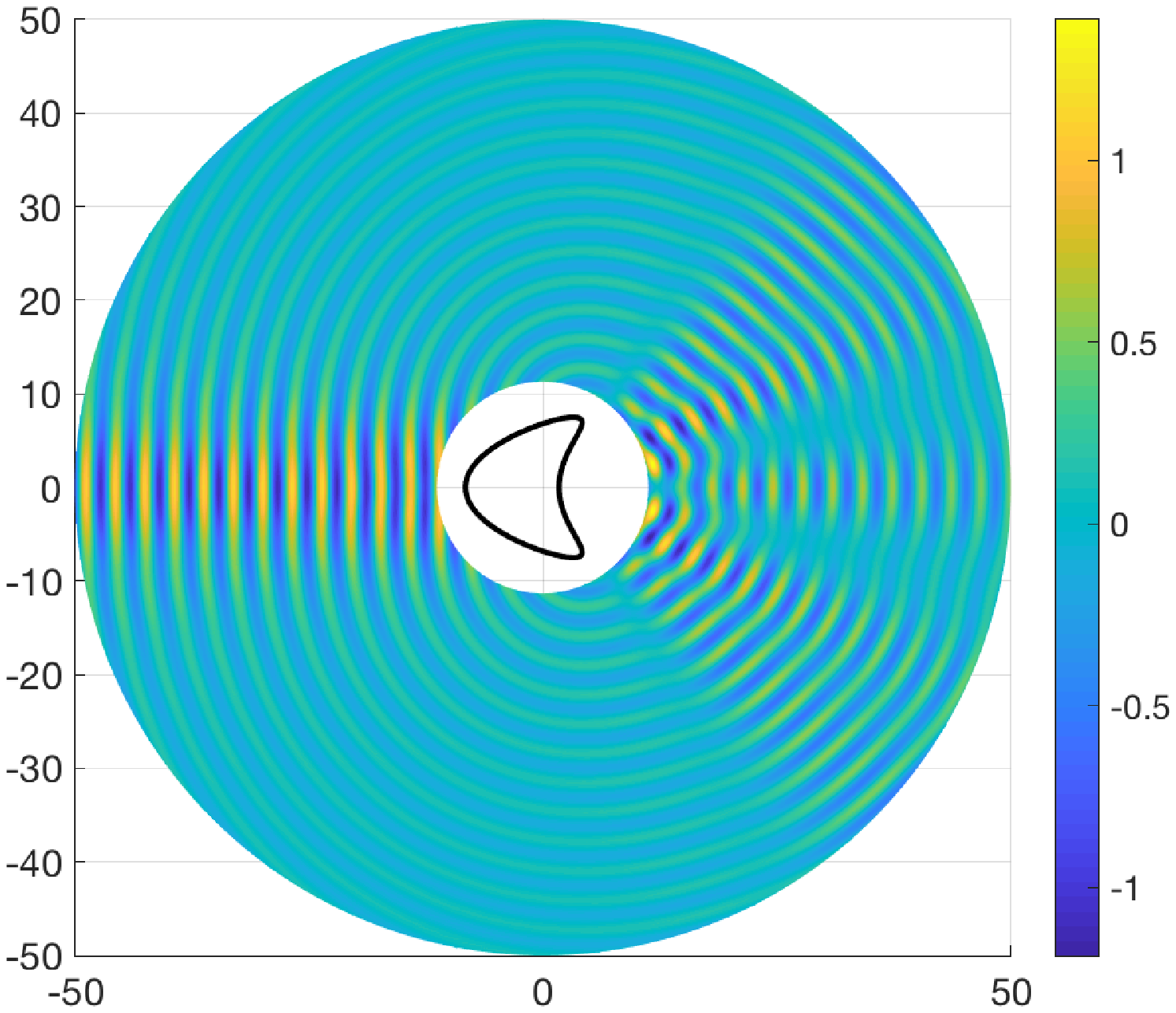} &
\includegraphics[trim = 65 25 65 25, clip,width=0.27\textwidth]{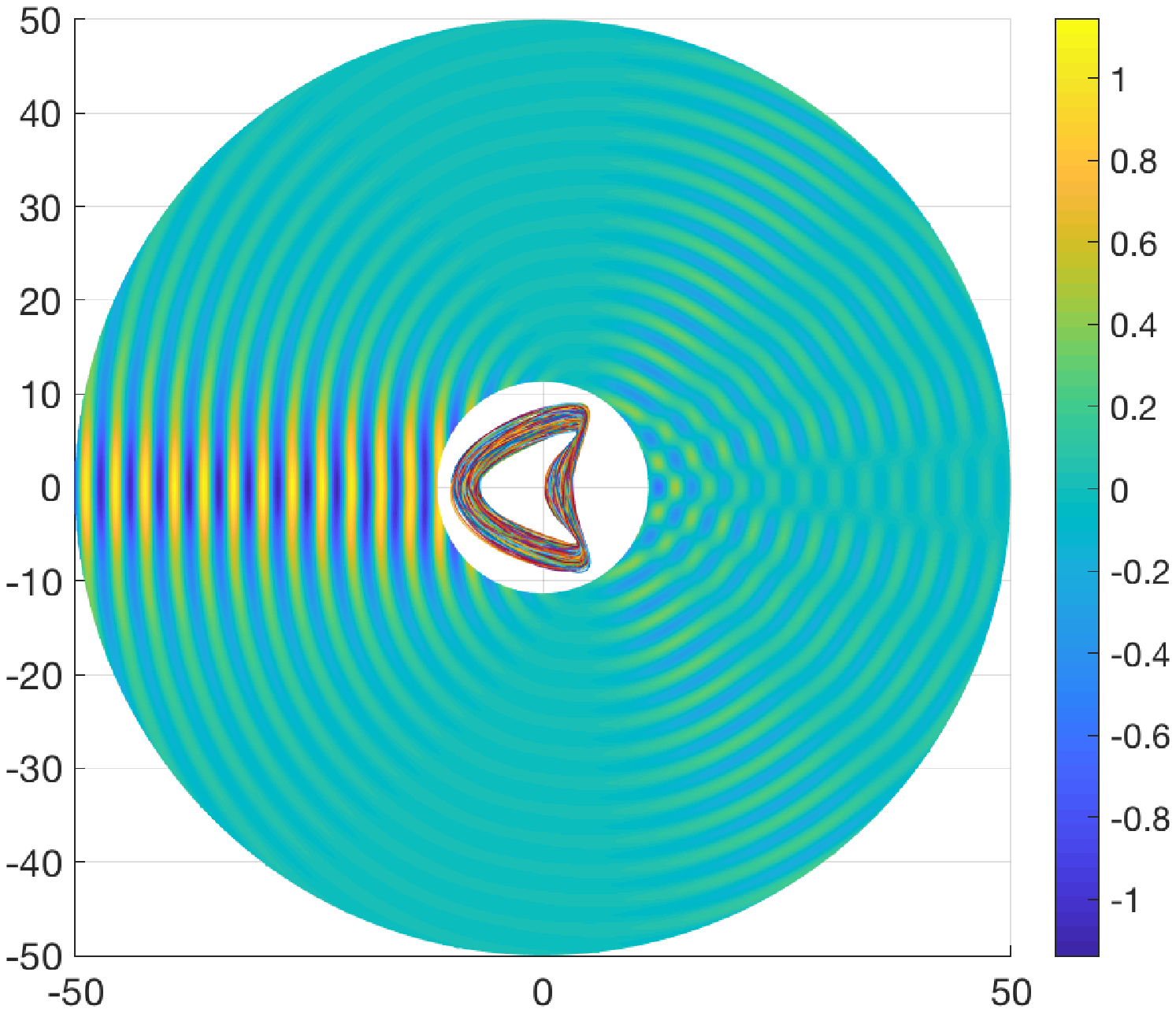}&
\includegraphics[trim = 65 25 65 25, clip,width=0.27\textwidth]{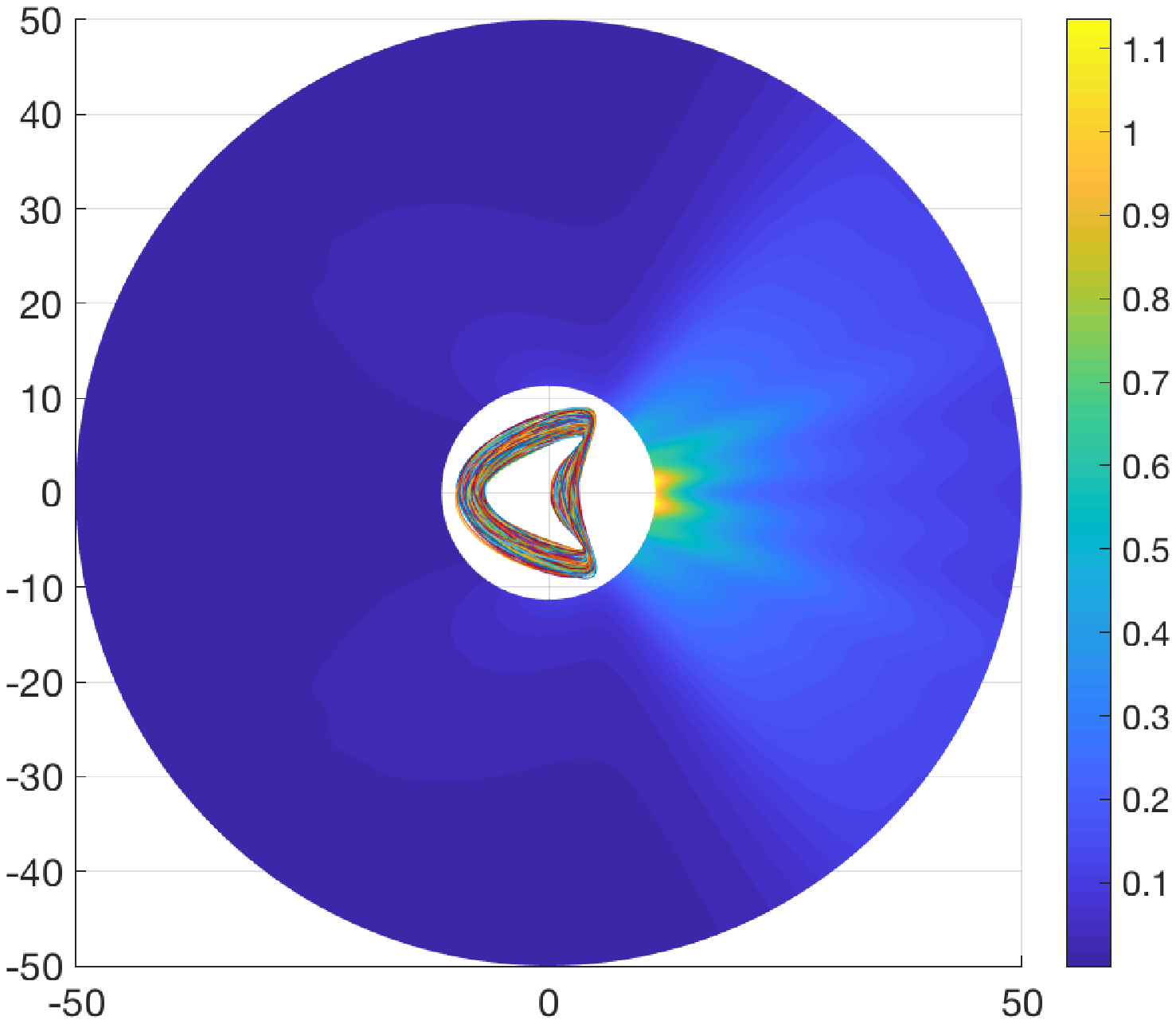}\\[0em]
\hline & & & \\[-1.1em]
\raisebox{3.4em}{\hspace*{-1em}\({\bf d}=\begin{bmatrix}-\nicefrac{1}{\sqrt{2}}\\\phantom{-}\nicefrac{1}{\sqrt{2}}\end{bmatrix}\)\hspace*{-0.5em}}&
\includegraphics[trim = 65 25 65 25, clip,width=0.27\textwidth]{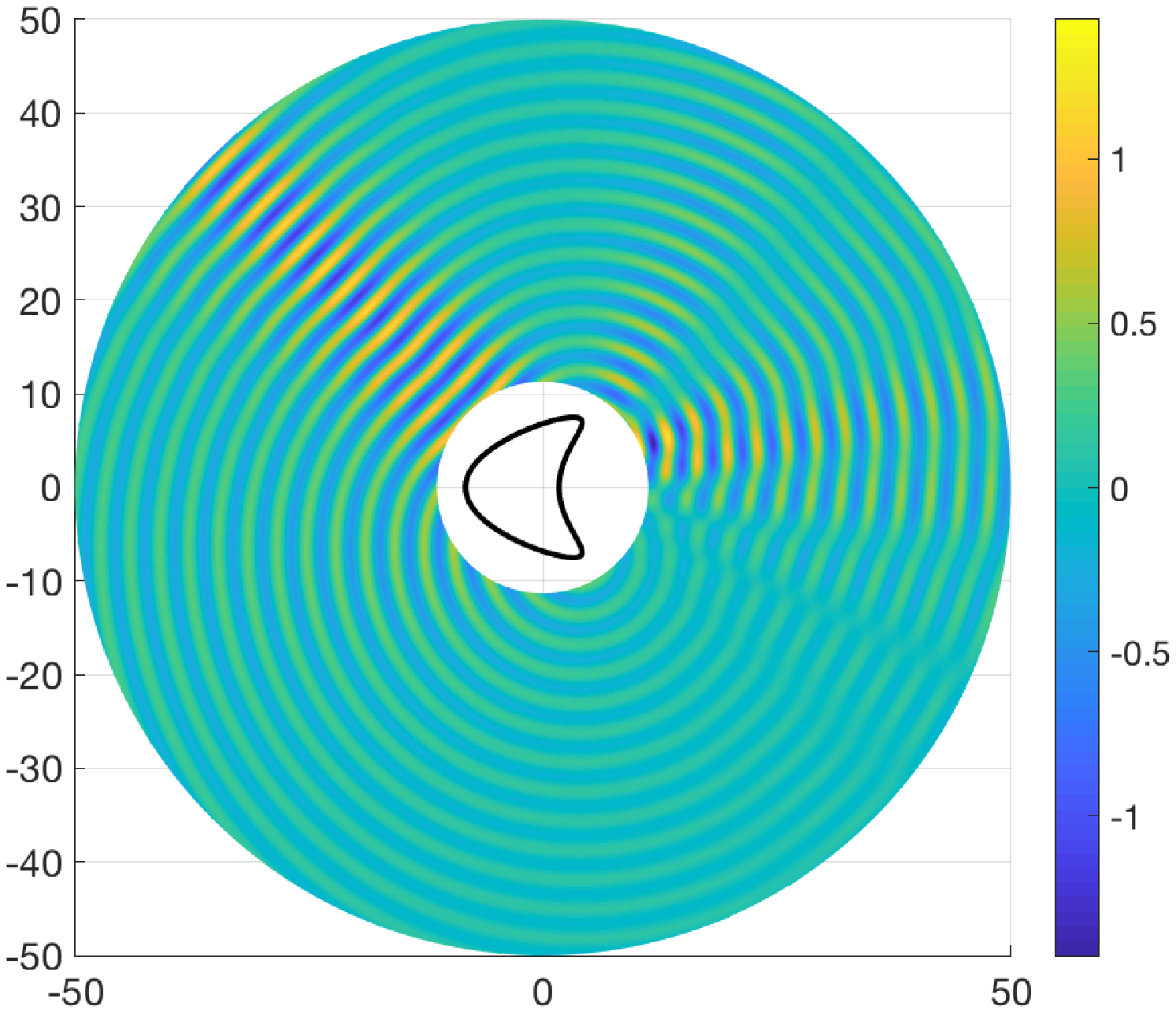} &
\includegraphics[trim = 65 25 65 25, clip,width=0.27\textwidth]{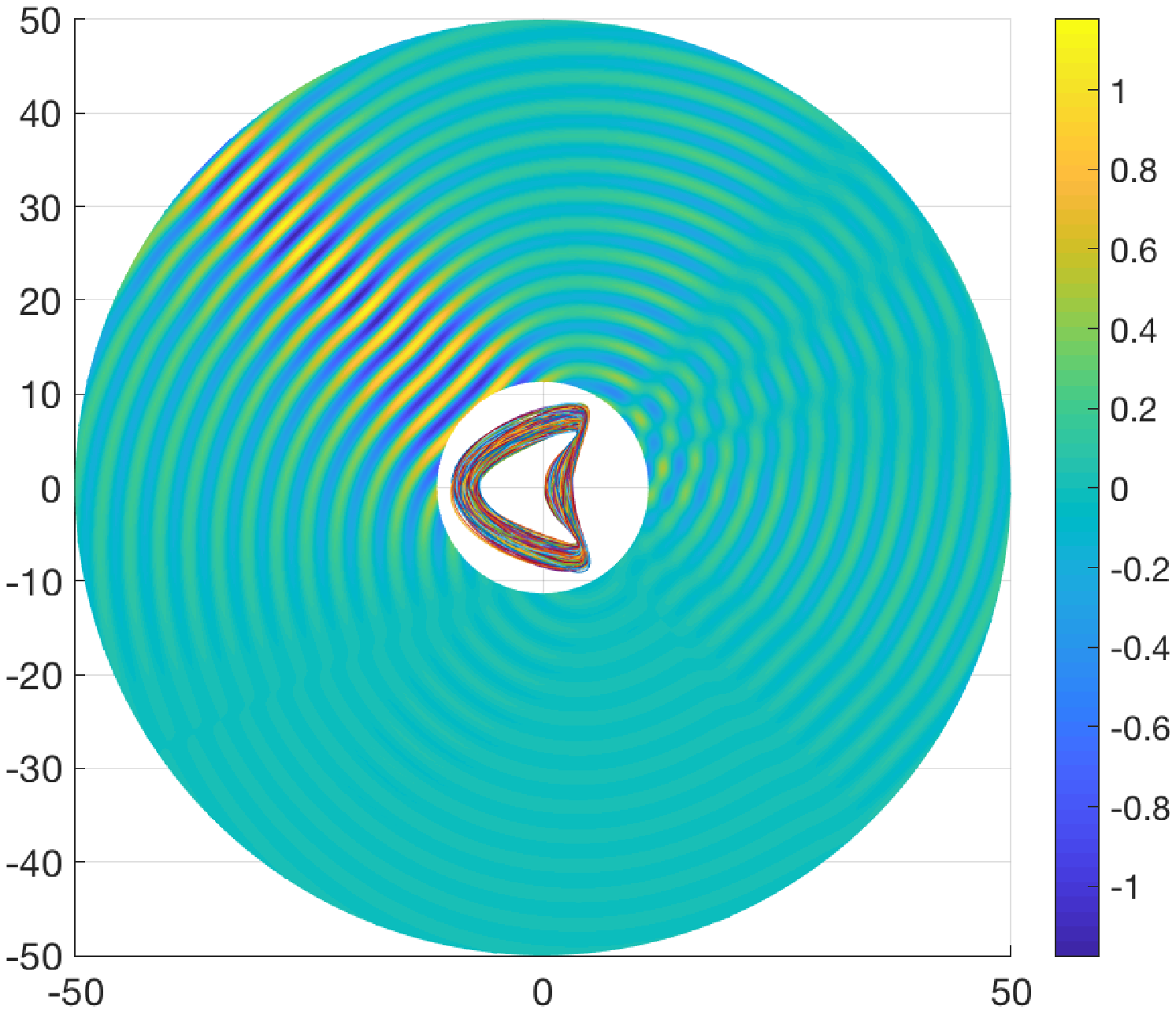}&
\includegraphics[trim = 65 25 65 25, clip,width=0.27\textwidth]{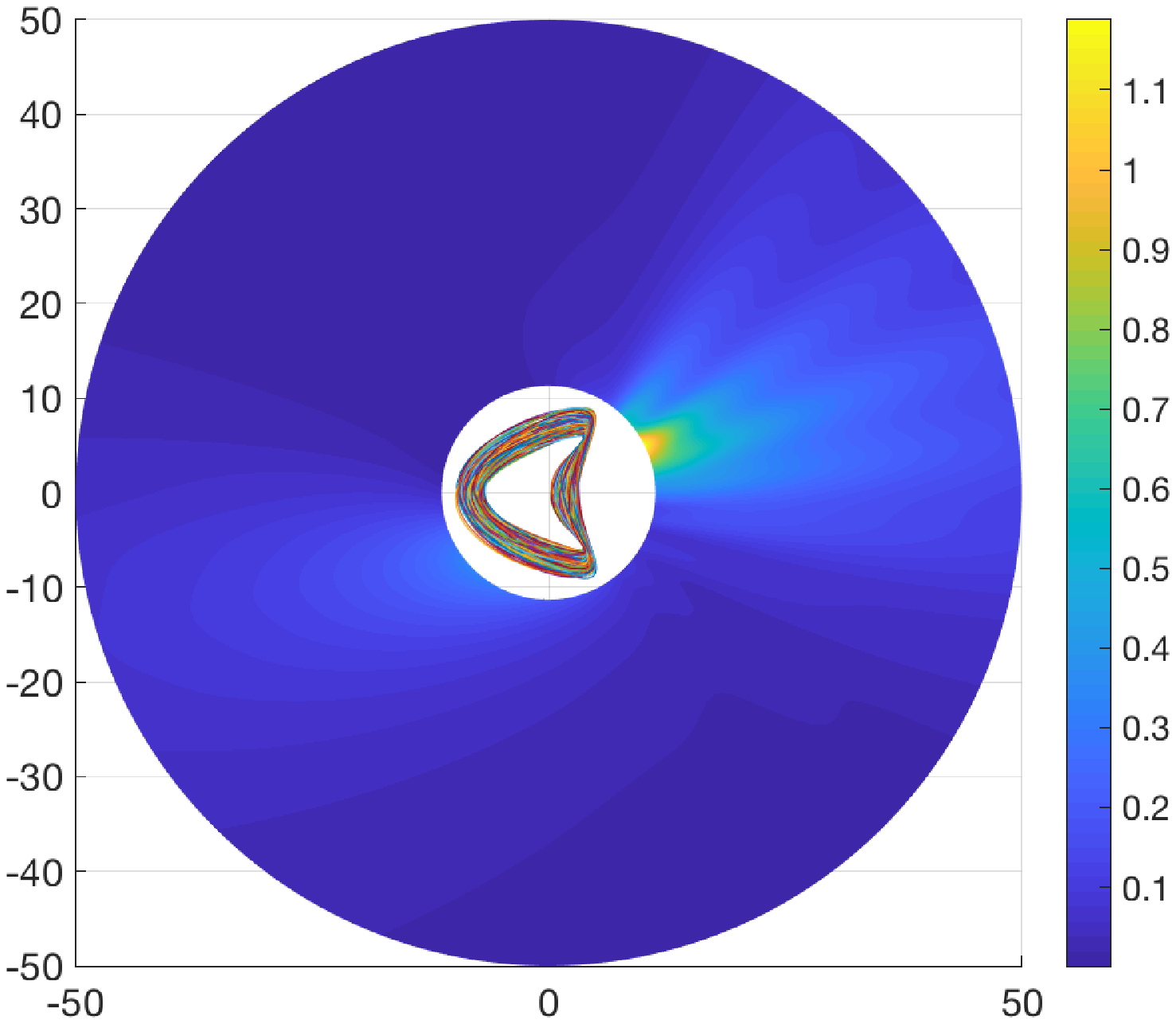}\\[0em] \hline & & & \\[-1.1em]
\raisebox{3.4em}{\hspace*{-1em}\({\bf d}=\begin{bmatrix}\phantom{-}0\\ -1\end{bmatrix}\)\hspace*{-0.5em}}&
\includegraphics[trim = 65 25 65 25, clip,width=0.27\textwidth]{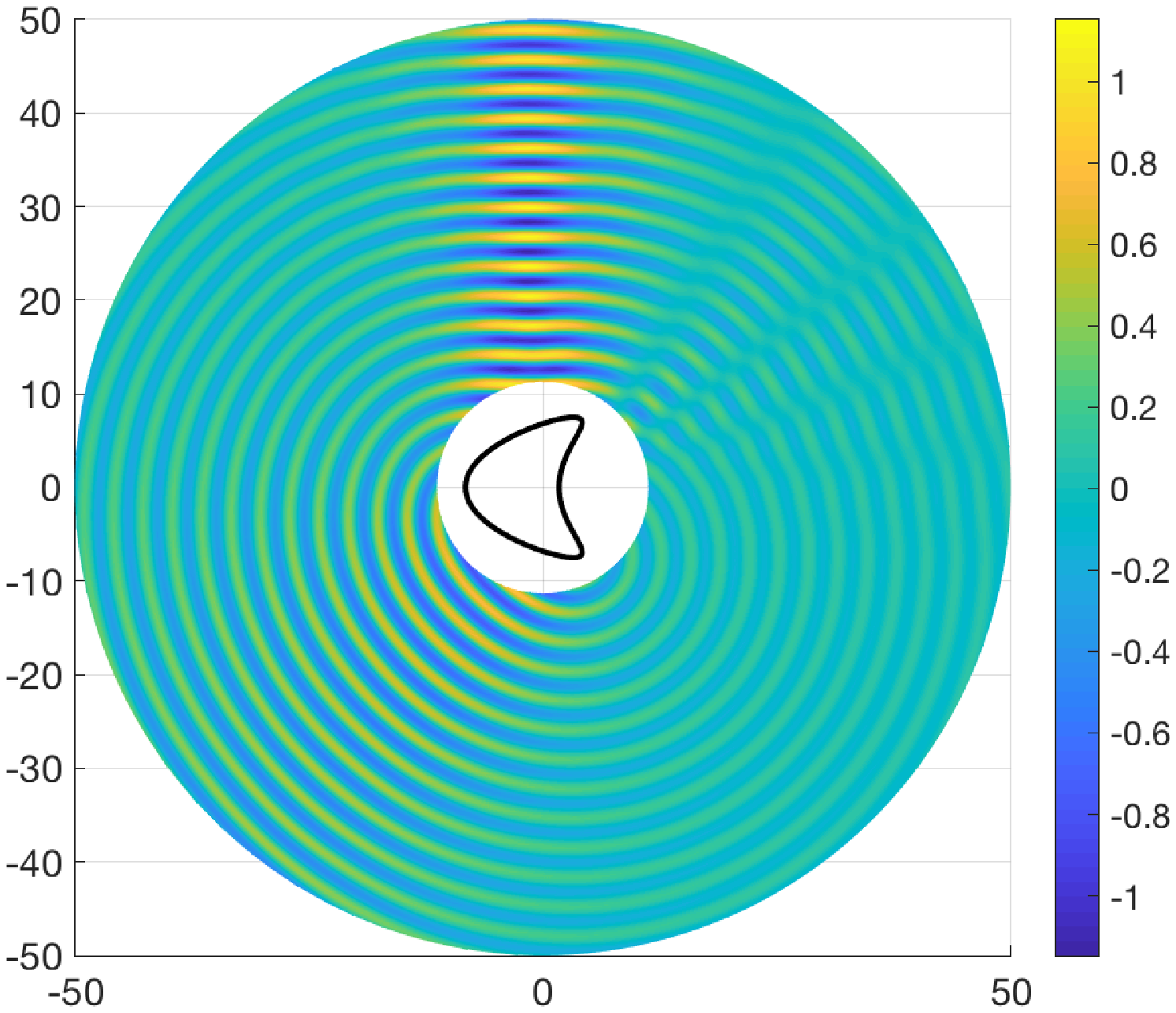} &
\includegraphics[trim = 65 25 65 25, clip,width=0.27\textwidth]{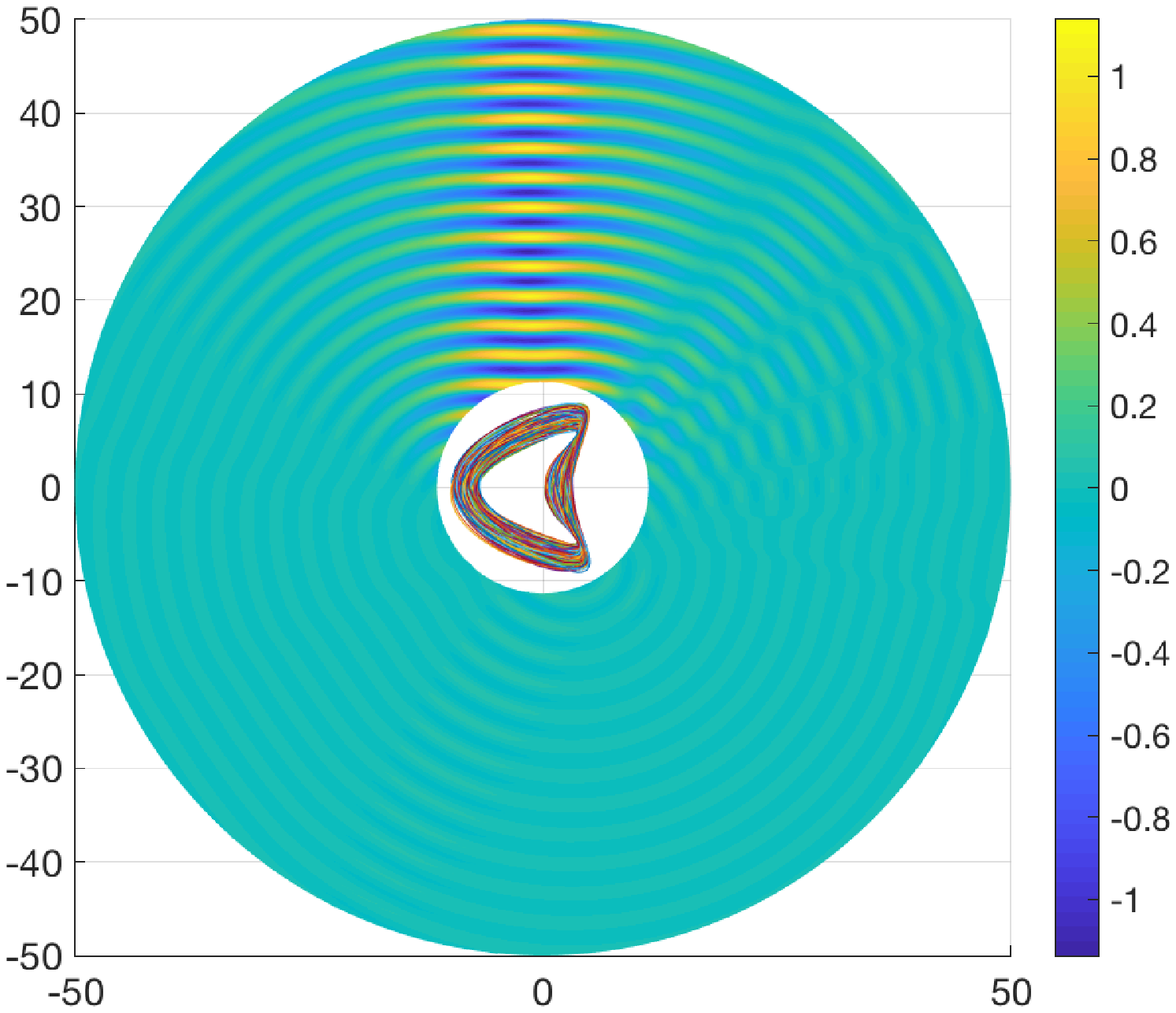}&
\includegraphics[trim = 65 25 65 25, clip,width=0.27\textwidth]{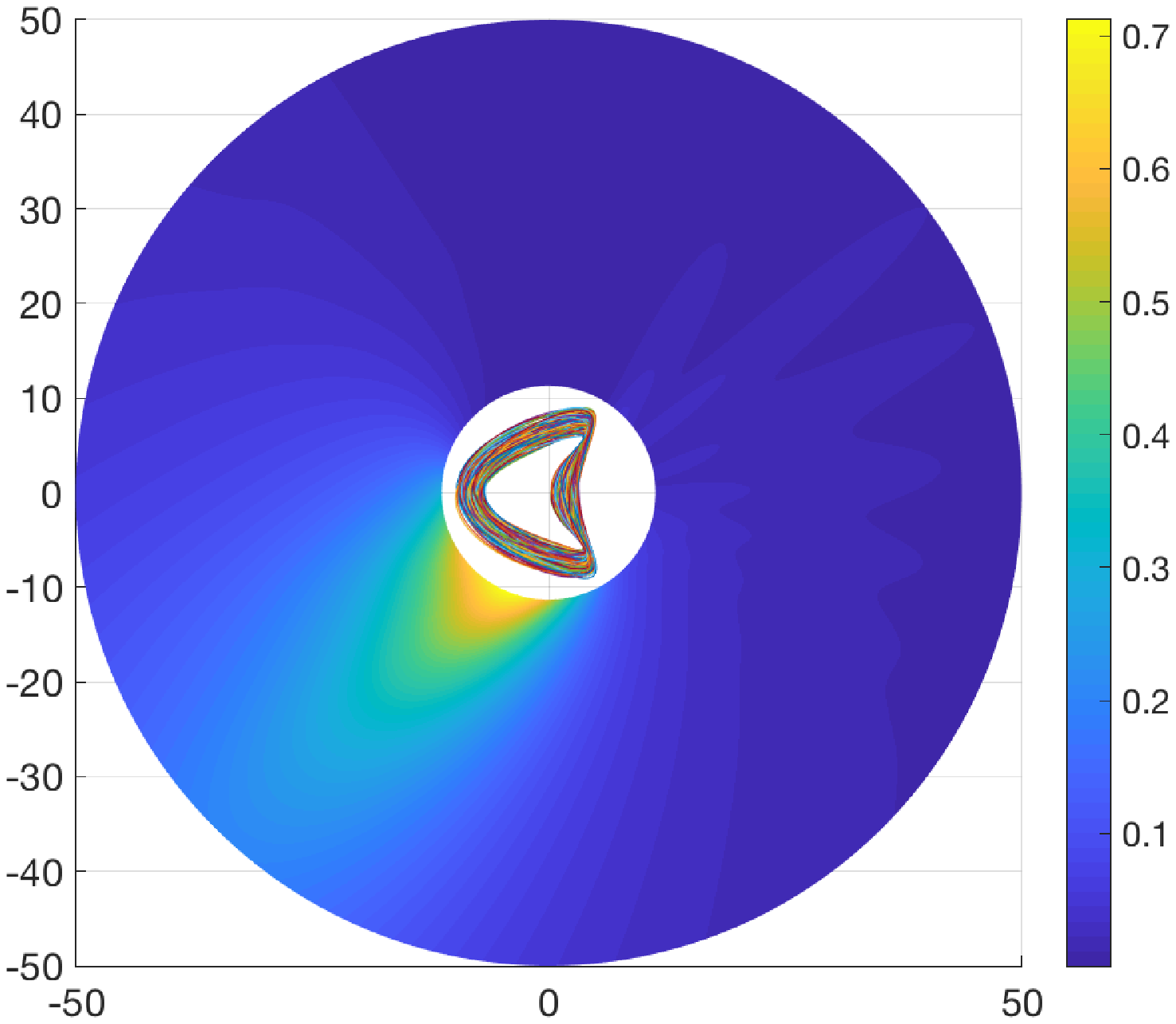}\\[0em] \hline & & & \\[-1.1em]
\raisebox{3.4em}{\hspace*{-1em}\({\bf d}=\begin{bmatrix}\nicefrac{1}{\sqrt{2}}\\\nicefrac{1}{\sqrt{2}}\end{bmatrix}\)\hspace*{-0.5em}}&
\includegraphics[trim = 65 25 65 25, clip,width=0.27\textwidth]{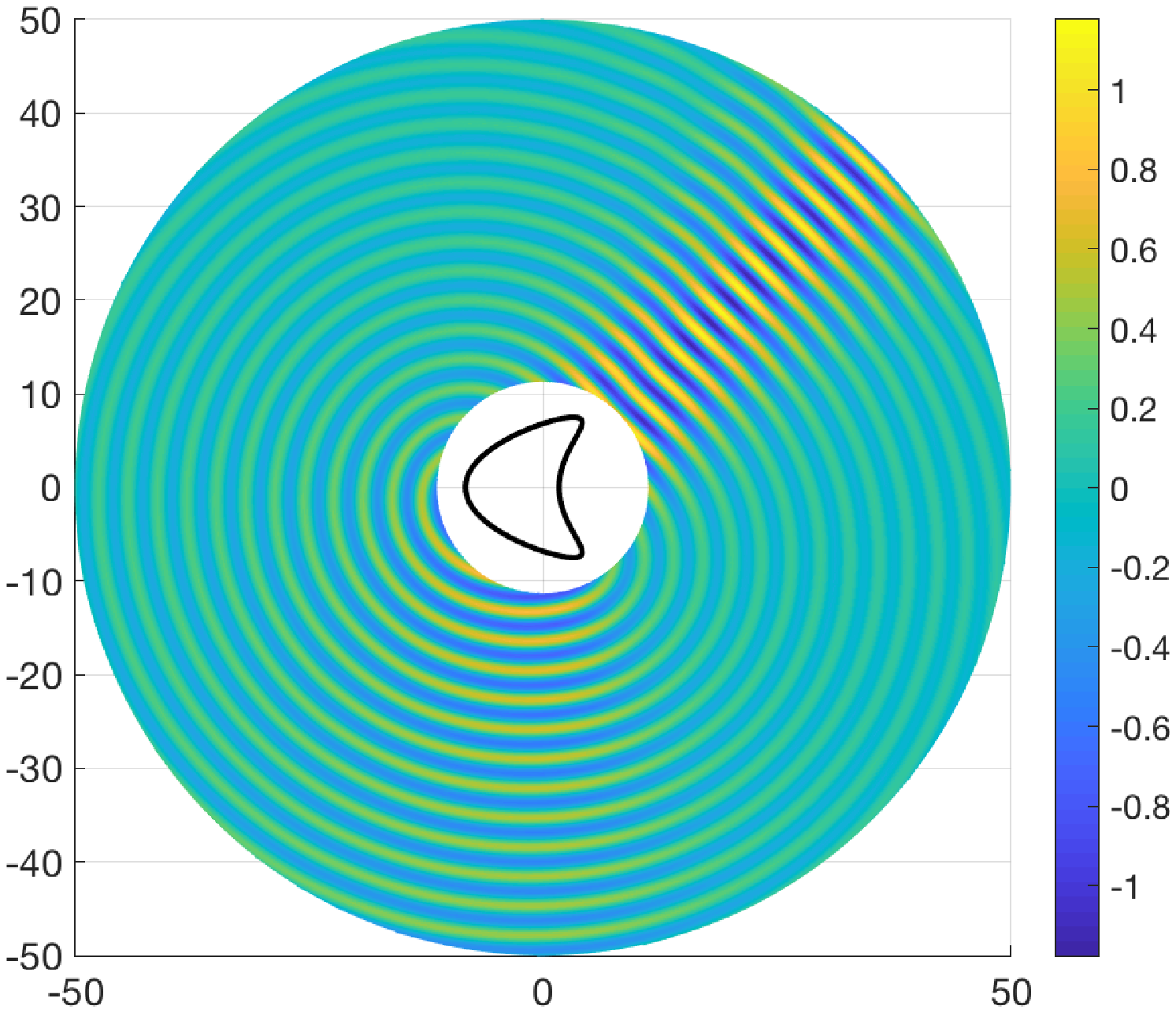} &
\includegraphics[trim = 65 25 65 25, clip,width=0.27\textwidth]{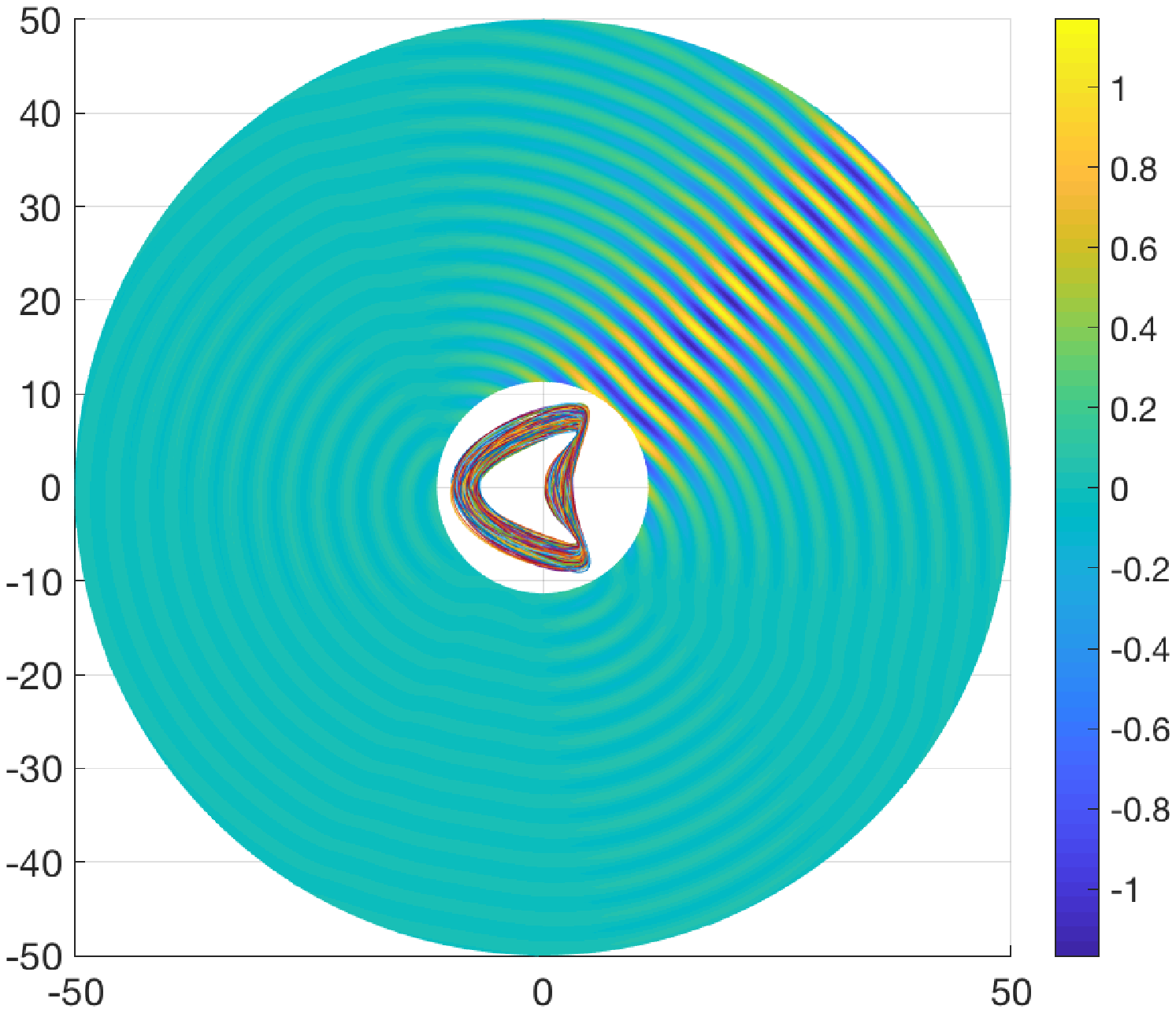}&
\includegraphics[trim = 65 25 65 25, clip,width=0.27\textwidth]{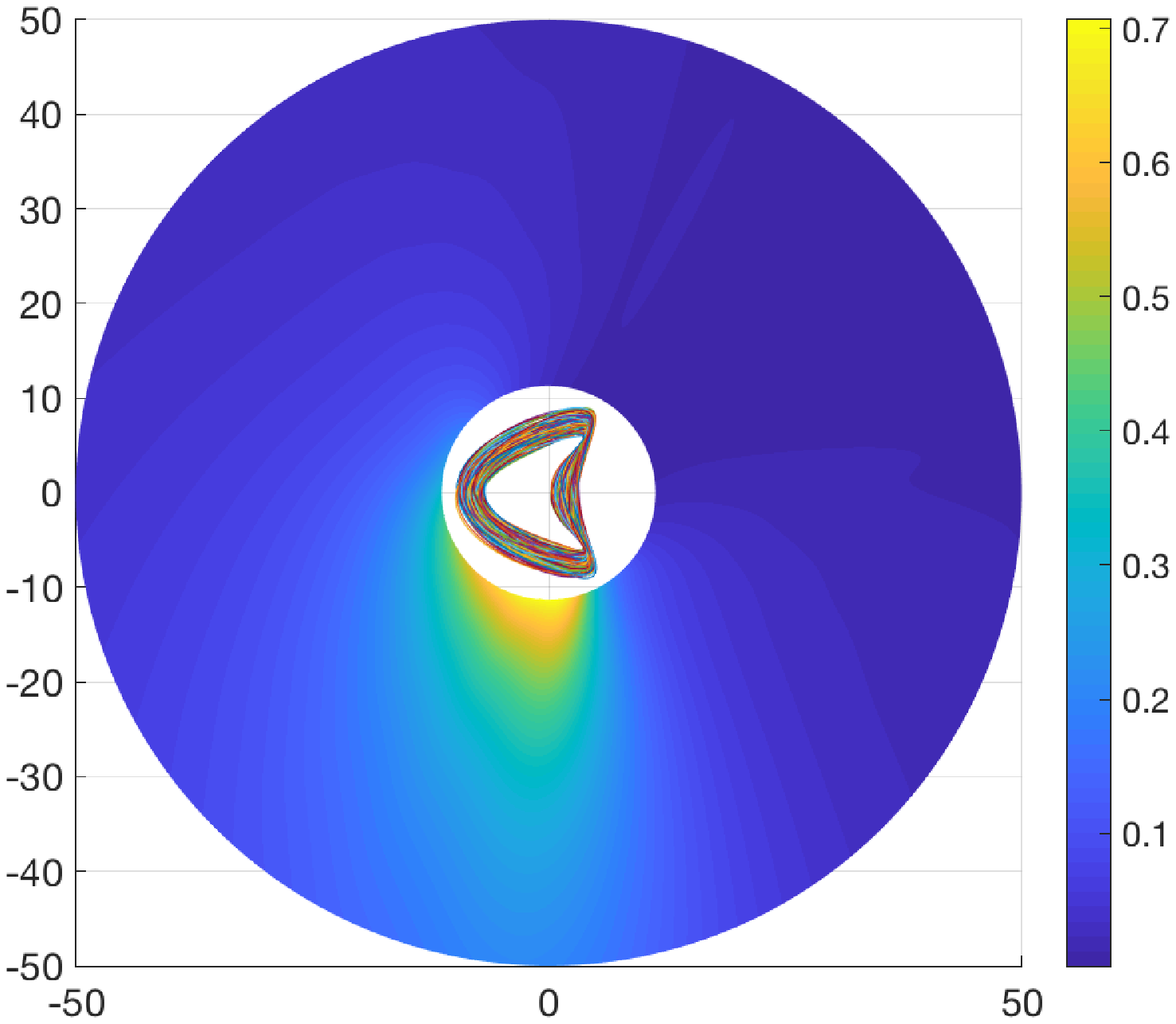}\\[0em]
\end{tabular}
\caption{\label{fig:direction}
The total field of the nominal scatterer (left), 
the expected total field for the random scatterer (middle) and
associated variance (right) for different directions of the incident wave.}
\end{center}
\end{figure}

In the previous paragraph, we have computed the second order 
statistics of the scattered wave for varying wavenumber, 
while the incident wave was always fixed to ${\bf d} = [1,0]^\top$. 
Now, we shall fix the wavenumber to $\kappa = 2$ and 
consider different directions of the incident wave. The results are 
depicted in Figure \ref{fig:direction}. Here, the first row
shows the total field for the nominal scatterer, the expected total 
wave for a random scatterer and its variance for 
an incident wave from the right in the first row,
for an incident wave from the bottom-right in the second row,
for an incident wave from the bottem in the third row and
finally for an incident wave from the bottom-left in the last row.
It turns out that the total field is mostly affected by
the perturbation of the scatterer in the direction of
the incident wave. In particular, we observe a high variance, where
the incoming wave hits the scatterer, while the variance
is nearly zero in the shadow region.
 
\subsection{Far-field pattern}
\begin{figure}[htb]
\begin{center}
\begin{minipage}{0.48\textwidth}
\begin{center}
\(\kappa=1\)\\[0.2em]
\includegraphics[clip,width=\textwidth]{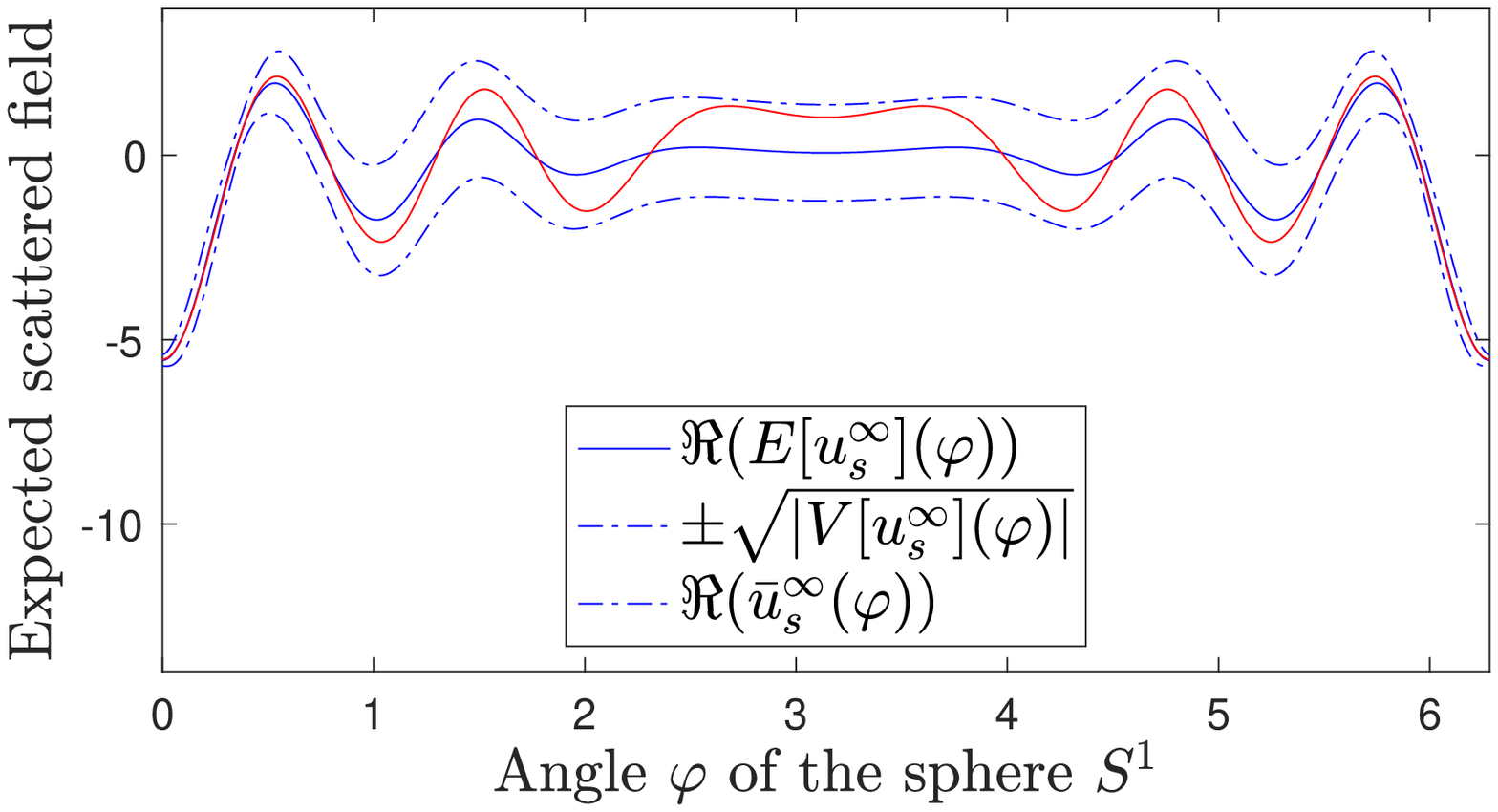}
\end{center}
\end{minipage}\hfill
\begin{minipage}{0.48\textwidth}
\begin{center}
\(\kappa=2\)\\[0.2em]
\includegraphics[clip,width=\textwidth]{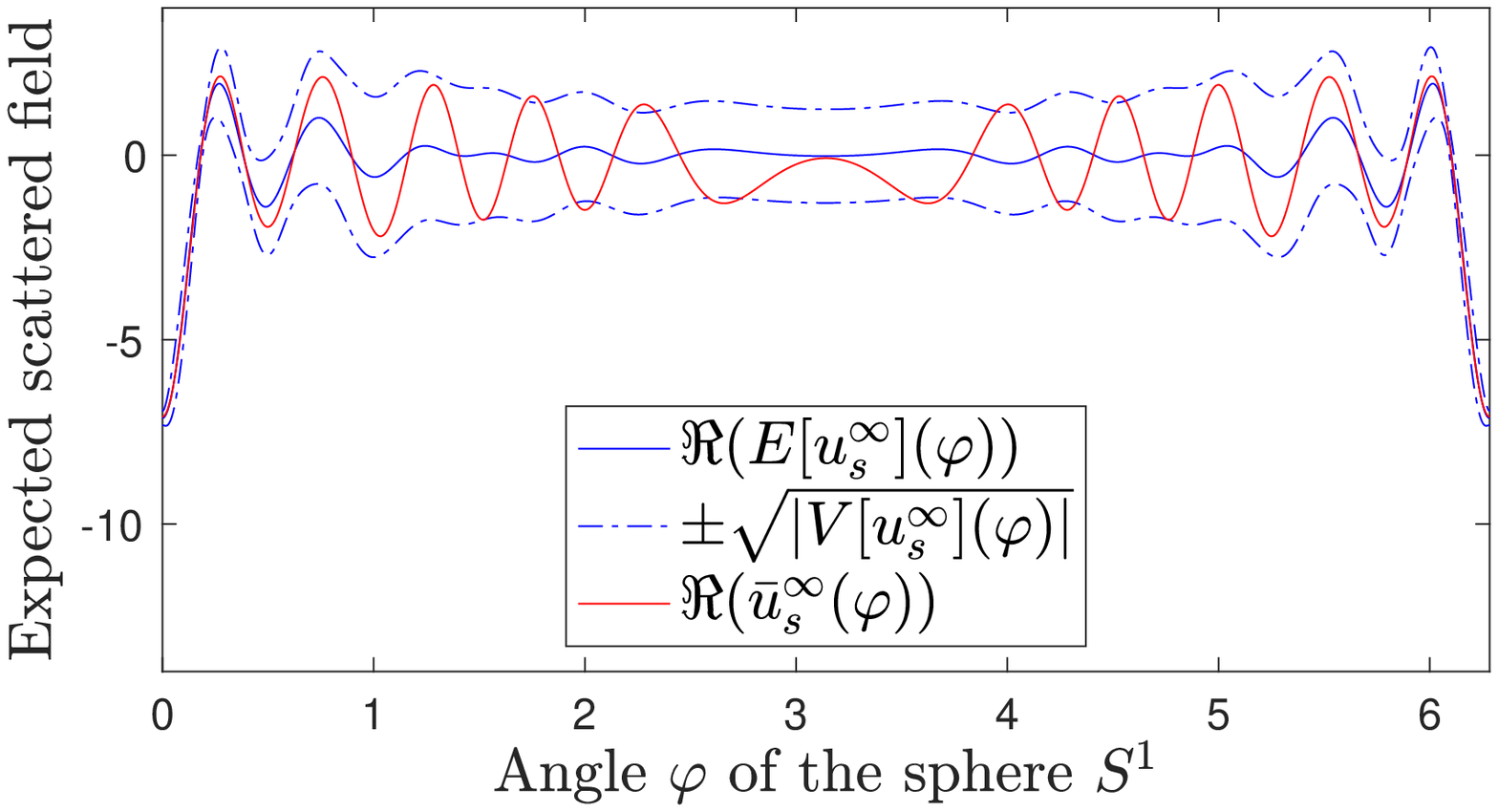}
\end{center}
\end{minipage}\\[1em]
\begin{minipage}{0.48\textwidth}
\begin{center}
\(\kappa=4\)\\[0.2em]
\includegraphics[clip,width=\textwidth]{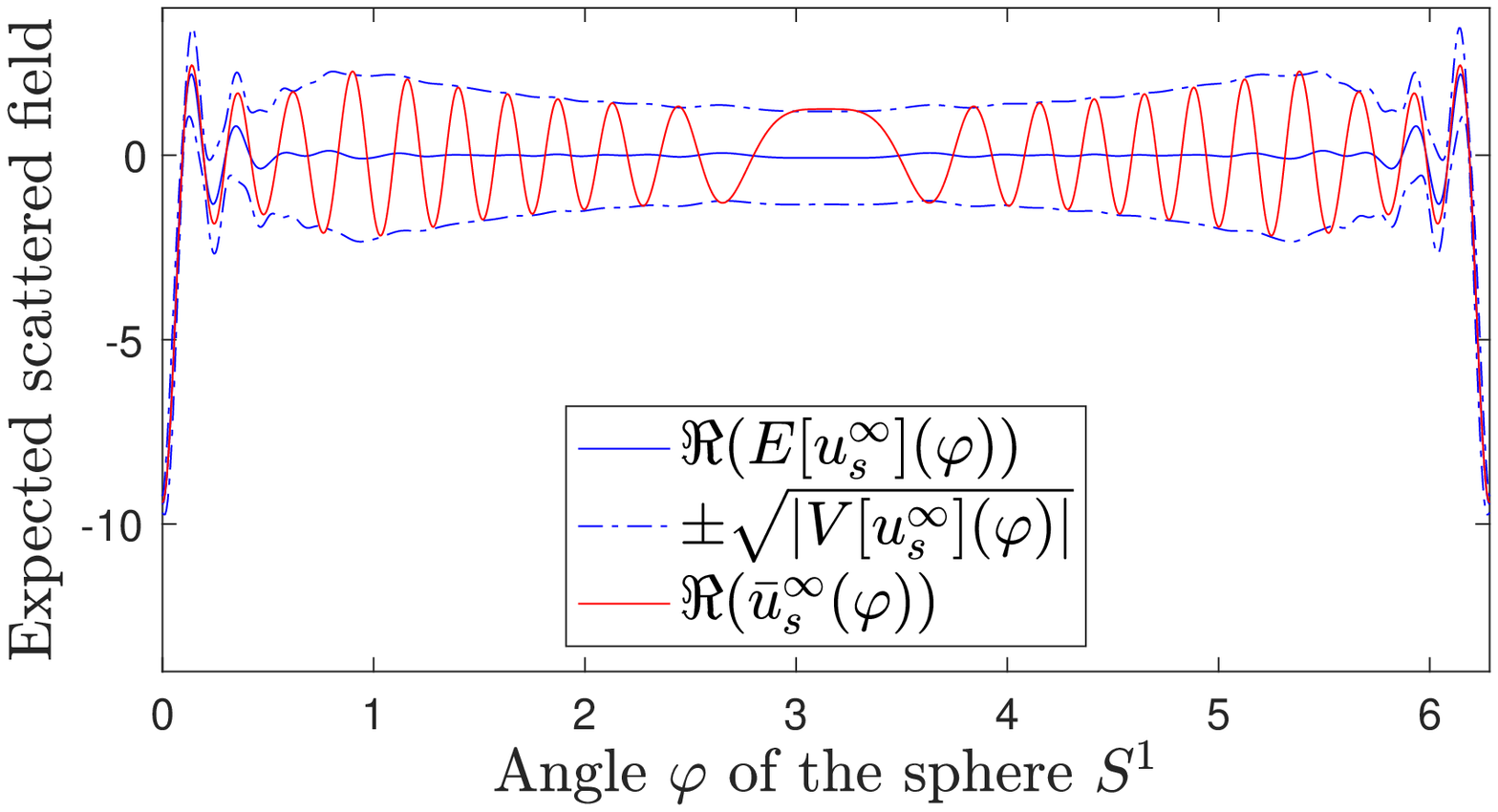}
\end{center}
\end{minipage}\hfill
\begin{minipage}{0.48\textwidth}
\begin{center}
\(\kappa=8\)\\[0.2em]
\includegraphics[clip,width=\textwidth]{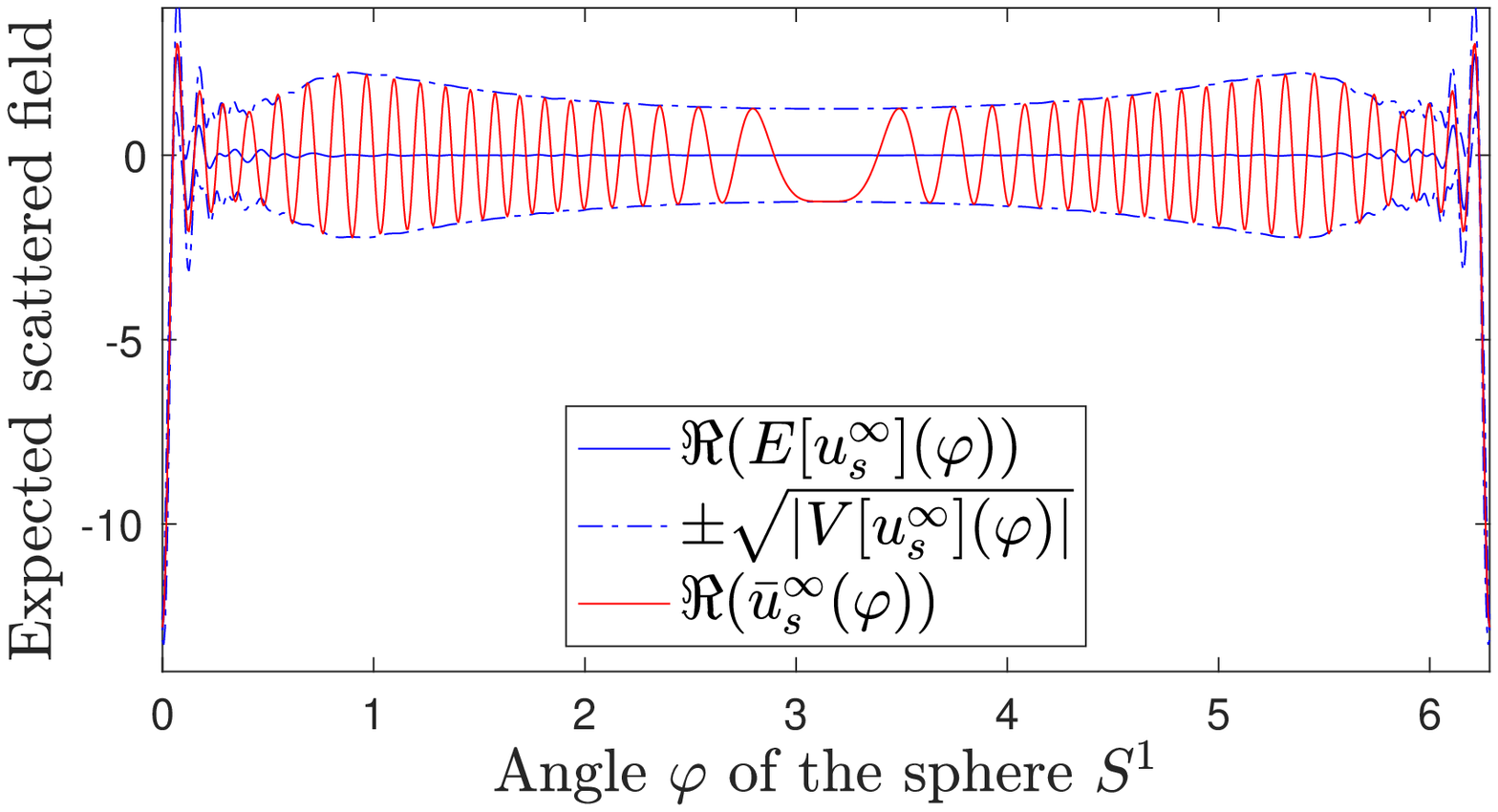}
\end{center}
\end{minipage}
\caption{\label{fig:far-field}The real part of the expected 
far-field (blue line) with standard deviation (blue dash-dotted line) and
the far-field of the unperturbed kite shape (red line) for
\(\kappa=1,2,4,8\).}
\end{center}
\end{figure}

We shall next consider the far-field pattern of the randomly 
perturbed kite-shape scatterer for the wavenumbers 
$\kappa=1,2,4,8$. The far-field has been evaluated 
in $n = 1000$ equidistant points on $\mathbb{S}^1$. 
The results of the computations are found in 
Figure~\ref{fig:far-field}, 
where we have depicted the expected far-field
(blue line) and the standard deviation of the far-field (dash-dotted 
line). We observe that the expected far-field pattern 
only oscillates in the shadow region. \red{More precisely,
on the left of the scatterer, i.e.\ $\pi/2< \varphi <3\pi/2$, 
the average perturbed far-field does 
not follow the oscillations in the far-field of the nominal geometry (red line). 
This is in contrast to the deep shadow range, obtained for 
$\varphi\approx 0$, where the average far-field and nominal 
far-field are rather close to each other.}
\section{Conclusion}\label{sec:conclusion}
\red{
In the present article, we have proposed an efficient method for
the computation of far-field statistics for acoustic scattering
in the case of random obstacles. We have employed the Karhunen-Lo\`eve 
expansion to parametrize the random scattering problem with 
respect to the infinite dimensional hypercube. Then, the parametric 
scattering problem has been reformulated by means of boundary 
integral equations. This approach directly leads to a reduction of 
the spatial dimensionality from \(d\) to \(d-1\) and, consequently, 
to a reduction of the computational cost.}

\red{
For the rapid computation of far-field statistics,
like the mean and the variance of the far-field pattern,
we have introduced an artificial boundary, 
which almost surely contains all realizations of the random scatterer.
Using this approach, all information of the random domain perturbation
is assessable from the scattered wave's Cauchy data.
In particular, expressions for the far-field's
or the scattered wave's expectation and variance can easily be derived.
In order to speed up the computation of the variance, which is
based on the evaluation of the covariance, we have suggested the
application of a low-rank approximation method, namely the
pivoted Cholesky decomposition.}

\red{
The presented approach can also be formulated for 
the scattering at sound-hard obstacles. Moreover, it can 
be extended to other boundary value problems for which
a Green's function is available. The observation that the 
solution's second order statistics is determined by the 
second order statistics of the Cauchy data on a deterministic 
interface holds even for arbitrary second order elliptic 
boundary value problems.}

\red{
Numerical results have been provided for \(d=2\). Here, the spatial
discretization has been performed by the exponentially convergent
Nystr\"om method, while the quadrature in the random parameter
is approximated by the Halton sequence.
The numerical approximation of the Cauchy data for \(d=3\) can
be facilitated with linear cost in terms of degrees of freedom by
a fast boundary element method, like e.g.\ fast multipole, see \cite{GR87}, 
or wavelets, see \cite{HS06}. Since the evaluation of the correlation
always involves the evaluation of the potential, a fast algorithm
is required here as well. This can also be realized by the use of 
the fast multipole method, see \cite{OSU10}.
}
\bibliographystyle{plain}
\bibliography{literature}
\end{document}